\documentclass[a4paper,10pt]{article}

\usepackage{bm}
\usepackage{pifont}

\usepackage{amsmath,amssymb,amsfonts,color}
\usepackage{graphicx,epsfig,epstopdf} 
\usepackage{array} 
\usepackage[usenames,dvipsnames]{xcolor}
\usepackage{appendix}

\usepackage{indentfirst} 
\usepackage{enumerate} 
\usepackage{placeins} 
\usepackage{booktabs,ctable,multirow,longtable} 
\usepackage{psfrag}
\usepackage{graphicx}

\definecolor{jade}{rgb}{0.0, 0.66, 0.42}
\definecolor{darkorange}{rgb}{1.0, 0.55, 0.0}
\definecolor{fuchsia}{rgb}{1.0, 0.0, 1.0}
\definecolor{brown(web)}{rgb}{0.65, 0.16, 0.16}
\definecolor{dukeblue}{rgb}{0.0, 0.0, 0.61}

\newtheorem{remark}{Remark}
\makeatother

\begin{document}

\title{On penalization in variational phase-field models of brittle fracture}

\author{Tymofiy Gerasimov$^{1,\star}$, Laura De Lorenzis$^1$}

\maketitle

\begin{center} \footnotesize
$^1$ Institute of Applied Mechanics, Technische Universit\"at Braunschweig, Pockelsstr.\ 3, Braunschweig 38106, Germany \\	
$^\star$ Correspondence: t.gerasimov@tu-braunschweig.de
\end{center}

\begin{abstract}
Irreversible evolution is one of the central concepts as well as implementation challenges of both the variational approach to fracture by Francfort and Marigo (1998) and its regularized counterpart by Bourdin, Francfort and Marigo (2000, 2007 and 2008), which is commonly referred to as a phase-field model of brittle fracture. Irreversibility of the crack phase-field imposed to prevent fracture healing leads to a constrained minimization problem, whose optimality condition is given by a variational {\em inequality}. In our study, the irreversibility is handled via penalization. Provided the penalty constant is well-tuned, the penalized formulation is a good approximation to the original one, with the advantage that the induced {\em equality}-based weak problem enables a much simpler algorithmic treatment. We propose an analytical procedure for deriving the {\em optimal} penalty constant, more precisely, its {\em lower bound}, which guarantees a sufficiently accurate enforcement of the crack phase-field irreversibility. Our main tool is the notion of the optimal phase-field profile, as well as the $\Gamma$-convergence result. It is shown that the explicit lower bound is a function of two formulation parameters (the fracture toughness and the regularization length scale) but is independent on the problem setup (geometry, boundary conditions etc.) and the formulation ingredients (degradation function, tension-compression split etc.). The optimally-penalized formulation is tested for two benchmark problems, including one with available analytical solution. We also compare our results with those obtained by using the alternative irreversibility technique based on the notion of the history field by Miehe et al.\ (2010). \\

\end{abstract}

\tableofcontents

\section{Introduction}
\label{Intro}
\begin{figure}[!ht]
\begin{center}
\includegraphics[width=1.0\textwidth]{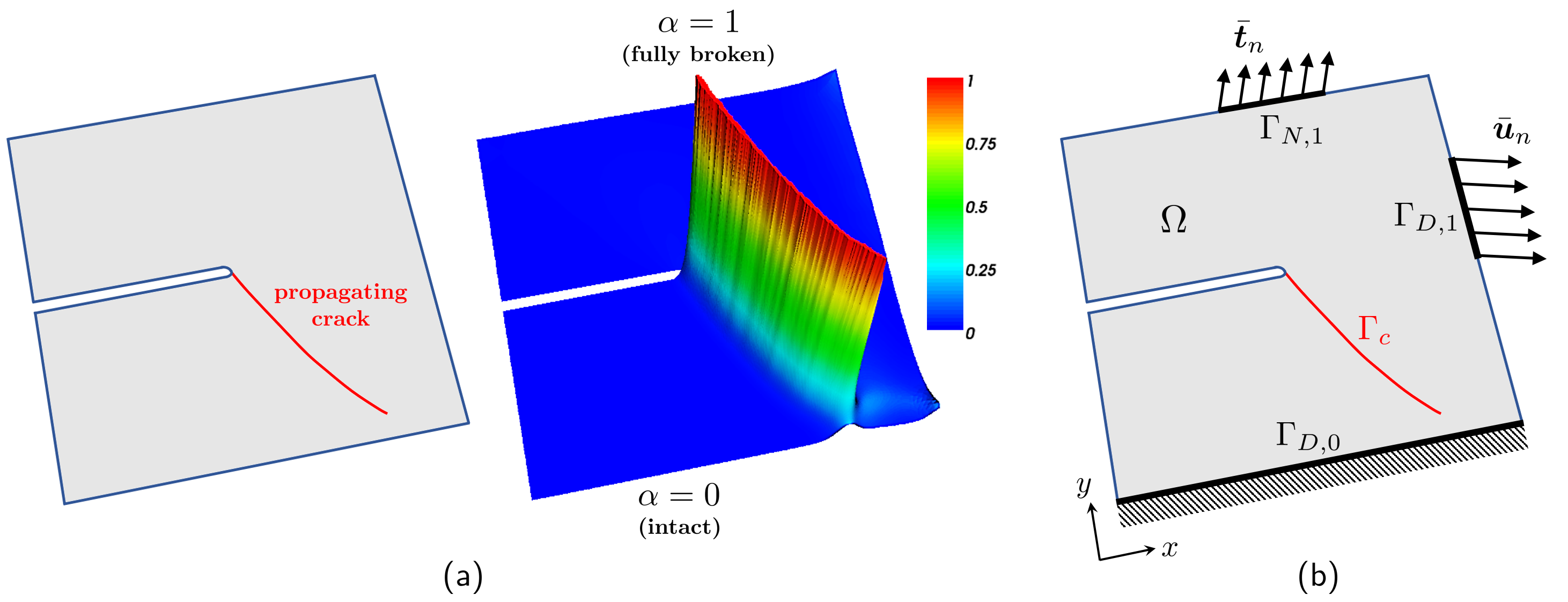}
\end{center}
\caption{(a) Phase-field description of fracture (sketchy) with $\alpha\in C(\Omega,[0,1])$ as the crack phase-field; (b) mechanical system setup used in section \ref{SecPF}.}
\label{Fig1}
\end{figure}

The phase-field framework for modeling systems with sharp interfaces consists in incorporating a continuous field variable -- the so-called order parameter -- which differentiates between multiple physical phases within a given system through a smooth transition. In the context of fracture, such an order parameter (termed the {\em crack phase-field}) describes the smooth transition between the fully broken and intact material phases, thus approximating the sharp crack discontinuity, as sketched in Figure \ref{Fig1}(a). The evolution of this field as a result of the external loading conditions models the fracture process.

The phase-field approach to brittle fracture dates back to the seminal work of Francfort and Marigo \cite{Francfort1998} on the {\em variational formulation of quasi-static brittle fracture} and to the related {\em regularized formulation} of Bourdin et al.\ \cite{Bourdin2000,Bourdin2007a,Bourdin2007b,Bourdin2008}. The former is the mathematical theory of quasi-static brittle fracture mechanics, which recasts Griffith's energy-based principle \cite{Griffith} as the minimisation problem of an energy functional. The latter presents an approximation, in the sense of $\Gamma$-convergence, of the energy functional and is designed to enable the numerical treatment.

The phase-field simulation of fracture processes holds a number of advantages over classical techniques with the discrete fracture description whose numerical implementation requires explicit (in the classical finite element method, FEM) or implicit (within the extended FEM) handling of the discontinuities. The most obvious one is the ability to track automatically a cracking process by the evolution of the smooth crack field on a fixed mesh. The possibility to avoid the tedious task of tracking complicated crack surfaces in 3d significantly simplifies the finite element implementation. The second advantage is the ability to simulate complicated processes, including crack initiation (also in the absence of a singularity), propagation, coalescence, branching and bifurcation without the need for additional ad-hoc criteria. With the formulation capability to also distinguish between fracture behavior in tension and compression, no supplementary contact problem has to be posed for preventing crack faces interpenetration. 

The currently available phase-field formulations of brittle fracture encompass static and dynamic models. We mention the papers by Del Piero et al.\ \cite{DelPiero2007}, Lancioni and Royer-Carfagni \cite{Lancioni2009},  Amor et al.\ \cite{Amor2009}, Freddi and Royer-Carfagni \cite{Freddi2009,Freddi2010}, Kuhn and M\"uller \cite{Kuhn2010}, Miehe et al.\ \cite{Miehe2010a,Miehe2010b}, Pham et al.\ \cite{Pham2011}, Borden \cite{Borden2012diss}, Borden et al.\ \cite{Borden2014}, Vignollet et al.\ \cite{Vignollet2014}, Mesgarnejad et al.\ \cite{Mesgarnejad2015}, Kuhn et al.\ \cite{Kuhn2015}, Ambati et al.\ \cite{Ambati2015review}, Strobl and Seelig \cite{Strobl2016}, Weinberg and Hesch \cite{Weinberg2017}, Tann\'e et al.\ \cite{{Tanne2018}}, Sargado et al.\ \cite{Sargado2018}, Gerasimov et al.\ \cite{Gerasimov2018}, where various formulations are developed and validated. Recently, the framework has been also extended to ductile (elasto-plastic) fracture \cite{Miehe2015b,Duda2015,Ambati2015,Alessi2015,Borden2016,Miehe2016,Alessi2018duct}, fracture in films \cite{Mesgarnejad2013,LeonBaldelli2014}, shells \cite{Shen2014,Ambati2016,Kiendl2016,Reinoso2017}, fracture under thermal loading \cite{Sicsic2013,Bourdin2014,Miehe2015a}, hydraulic fracture \cite{Bourdin2012hf,Wheeler2014,Mikelic2015a,Mikelic2015b,Mikelic2015c,Wilson2016}, fracture in porous media \cite{Miehe2015c,TaoWu2016,Cajuhi2017}, anisotropic fracture \cite{Li2015,Teichtmeister2017,Zhang2017,Nguyen2017a,Bleyer2018,Li2018}, fracture in laminates \cite{Alessi2018lam}, to name a few.

The finite element treatment of the formulation is, however, known to be computationally demanding. Two main reasons are the following: 
\begin{itemize}
	\item[(i)] the governing energy functional is non-convex and, in general, strongly non-linear with respect to both arguments (the displacement and the phase field). The staggered (also termed partitioned, or alternate minimization) solution approach based on decoupling of the strongly non-linear weak formulation into a system and then iterating between the equations is commonly used \cite{Bourdin2000,Bourdin2007a,Bourdin2007b,Bourdin2008,Amor2009,Miehe2010a,Miehe2010b,Pham2011,Borden2014,Mesgarnejad2015,Ambati2015review}. 
	The staggered scheme is robust, but typically has a very slow convergence behavior of the iterative solution process, see e.g.\ \cite{Ambati2015review,Gerasimov2016,Farrell2017}. Alternatively, the so-called monolithic approach which treats arguments simultaneously manifests major iterative convergence issues of the Newton-Raphson procedure due to non-convexity \cite{Gerasimov2016,Heister2015,Wick2016,Wick2017};
	\item[(ii)] the need to resolve the small length scale inherent to the diffusive crack approximation calls for extremely fine meshes, at least locally in the crack phase-field transition zone. Modeling a failure process whose final pattern is not known in advance precludes the construction of a suitably pre-refined mesh, thus forcing to compute on fixed uniform meshes (unless adaptivity is introduced). In this case, the computational cost is very high.
\end{itemize}
Already in the seminal paper by Bourdin et al.\ \cite{Bourdin2000} and later in \cite{Mesgarnejad2015} parallel computing has been advocated for the staggered solution scheme combined with uniformly fine meshes. However, some new results by Gerasimov and De Lorenzis \cite{Gerasimov2016}, Heister et al.\ \cite{Heister2015}, and Wick \cite{Wick2016,Wick2017} on the monolithic scheme, and by Farrell and Maurini \cite{Farrell2017} on over-relaxed accelerated staggered schemes hold a promise that efficient algorithmic handling of (i) is feasible. Furthermore, recent findings by Burke et al.\ \cite{Burke2010a,Burke2010b,Burke2013}, Artina et al.\ \cite{Artina2014,Artina2015} on error-controlled adaptive mesh refinement strategies, as well as by Heister et al.\ \cite{Heister2015} and Klinsmann et al.\ \cite{Klinsman2015} on physics-motivated procedures for mesh adaptivity provide a basis to efficiently tackle (ii) as well.

In this paper, we dissect another interesting modeling and computational ingredient of the formulation, namely, 
\begin{itemize}
	\item[(iii)] the {\em irreversibility} of the crack phase-field, i.e.\ the condition that prevents crack healing. It is given by the constraint $\dot{\alpha}\geq0$ in $\Omega$, with $\Omega$ being a computational domain. Using a backward difference quotient, this turns into $\alpha\geq\alpha_{n-1}$ in $\Omega$, with $n\geq1$ representing the pseudo-time or loading step in the incremental variational formulation.
\end{itemize}
Due to (iii), the formulation is a constrained minimization problem whose optimality condition is a {\em variational inequality} \cite{Amor2009,Pham2011,Mesgarnejad2015}, thus requiring special solution algorithms.  
Several options of enforcing the condition $\alpha\geq\alpha_{n-1}$ that lead to a simpler equality-based formulation are found in the literature. We classify these as relaxed, penalized and implicit ones. We first mention the works of Bourdin et al.\ \cite{Bourdin2000}--\cite{Bourdin2008}, where the irreversibility of $\alpha$ is enforced only on the so-called 'crack-set', that is, at the points of $\Omega$ where $\alpha=1$ or close to $1$. With this technique, irreversibility of only a fully developed crack is modeled and, therefore, it is viewed as a relaxed treatment of $\alpha\geq\alpha_{n-1}$. Bourdin et al.'s idea is adopted in \cite{DelPiero2007,Lancioni2009,Burke2010a} directly, in \cite{Burke2010b,Burke2013} with a modified notion of the crack-set and in \cite{Artina2014,Artina2015,Gerasimov2016} using penalization. With regard to the second option, we recall the augmented-Lagrangian method in \cite{Wheeler2014,Wick2016,Wick2017}. It yields the equality-based formulation equivalent to the original constrained minimization problem, but the presence of extra variables makes the computational effort high. Finally, the implicit enforcement of $\alpha\geq\alpha_{n-1}$ using a history field was proposed by Miehe et al.\ in \cite{Miehe2010b}. This method has been adopted in a major amount of works on the topic. The approach is, however, no longer of variational nature and its equivalence to the original problem cannot be proven. Also, in this case, only the staggered solution scheme can be employed.

In this manuscript we address the irreversibility constraint $\alpha\geq\alpha_{n-1}$ via simple penalization, with the advantages that the obtained equality-based formulation is equivalent to the original variational inequality problem (provided the penalty parameter is 'well chosen'), and that it can be considered with both staggered and monolithic solution schemes.

The paper is organized as follows. In section \ref{SecPF}, we outline the main concepts of phase-field modeling of brittle fracture and the formulation used in the present paper. The irreversibility constraint and its various implementation options are discussed in detail, including the penalized version in our focus. Section \ref{MainRes} is the main part of the manuscript. Therein, we start with recalling the notion of the {\em optimal} phase-field profile for a fully developed crack, as well as the related $\Gamma$-convergence result. The extension of this theory to our penalized formulation enables us to devise the analytical procedure for the 'reasonable' choice of a {\em lower bound} for the penalty parameter that guarantees a sufficiently accurate enforcement of the crack phase-field irreversibility. In section \ref{SecNum} we  present two numerical experiments that verify and illustrate our findings.

\section{Phase-field approach to brittle fracture}
\label{SecPF}
In this section, we consider a mechanical system undergoing a brittle fracture process modeled with the phase-field formulation and briefly recall the ingredients of the formulation and the algorithmic aspects of the solution.

\subsection{Governing energy functional}
Let $\Omega\subset\mathbb{R}^d$, $d=2,3$ be an open and bounded domain representing the configuration of a $d$-dimensional linear elastic body, and let $\Gamma_{D,0},\Gamma_{D,1}$ and $\Gamma_{N,1}$ be the (non-overlapping) portions of the boundary $\partial\Omega$ of $\Omega$ on which homogeneous Dirichlet, non-homogeneous Dirichlet and Neumann boundary conditions are prescribed, respectively. The body is assumed to be linearly elastic and isotropic, with the elastic strain energy density function given by $\Psi(\bm\varepsilon):=\frac{1}{2}\bm\varepsilon:\mathbb{C}:\bm\varepsilon=\frac{1}{2}\lambda\mathrm{tr}^2(\bm\varepsilon)+\mu\mathrm{tr}(\bm\varepsilon\cdot\bm\varepsilon)$, where, in turn, $\bm\varepsilon$ is the second-order infinitesimal strain tensor, $\mathbb{C}$ is the fourth-order elasticity tensor, and $\lambda$ and $\mu$ are the Lam\'{e} constants. Also, let $G_c$ be the material fracture toughness of the body. A quasi-static loading process with the discrete pseudo-time step parameter $n=1,2,...$, such that the displacement $\bar{\bm u}_n$ and traction $\bar{\bm t}_n$ loading data are prescribed on the corresponding parts of the boundary is considered. Finally, let $\Gamma_c\subset\Omega$ be the crack surface that is evolving during the process, see Figure \ref{Fig1}(b). 

For the mechanical system at hand, the variational approach to brittle fracture in \cite{Francfort1998} relies on the energy functional
\begin{equation}
\mathcal{E}(\bm u,\Gamma_c) = 
\underbrace{ \int_{\Omega\setminus\Gamma_c} 
\Psi(\bm\varepsilon(\bm u)) \, \mathrm{d}{\bf x}}_{\displaystyle E_\mathrm{el.}(\bm u,\Gamma_c)}
+ \underbrace{ G_c\int_{\Gamma_c} \mathrm{d}s}_{\displaystyle E_S(\Gamma_c)}
-\int_{\Gamma_{N,1} } \bar{\bm t}_n \cdot \bm u \, \mathrm{d}s,
\label{VAF}
\end{equation}
with $\bm u:\Omega\backslash\Gamma_c\rightarrow\mathbb{R}^d$ as the displacement field and $\Gamma_c$ as the crack set, and the related minimization problem at each $n\geq 1$. In (\ref{VAF}), the functionals termed $E_\mathrm{el.}$ and $E_S$ represent the elastic energy stored in the body and the fracture surface energy dissipated within the fracture process. The latter rigorously reads $E_S(\Gamma_c)=G_c\mathcal{S}^{d-1}(\Gamma_c)$ with $\mathcal{S}^p$ as the so-called $p$-dimensional Hausdorff measure of the crack set $\Gamma_c$.  In simple terms, $\mathcal{S}^1(\Gamma_c)$ and $\mathcal{S}^2(\Gamma_c)$ represent the length and the surface area of $\Gamma_c$ when $d=2$ and $3$, respectively. In the following, we use for $\mathcal{S}^{d-1}(\Gamma_c)$ the simpler notation $\left|\Gamma_c\right|$.

The regularization of (\ref{VAF}) {\em \'a la} Bourdin-Francfort-Marigo \cite{Bourdin2000,Bourdin2007a,Bourdin2007b,Bourdin2008}, which is the basis for a variety of fracture phase-field formulations, reads as follows:
\begin{equation}
\mathcal{E}(\bm u,\alpha) = 
\underbrace{ \int_\Omega
{\sf g}(\alpha)\Psi(\bm\varepsilon(\bm u)) \, \mathrm{d}{\bf x}}_{\displaystyle E_\mathrm{el.}(\bm u,\alpha)} 
+ \underbrace{ \frac{G_c}{c_{\sf w}} \int_\Omega 
\left(\frac{{\sf w}(\alpha)}{\ell}+\ell|\nabla\alpha|^2\right)
\mathrm{d}{\bf x}}_{\displaystyle E_S(\alpha)}
-\int_{\Gamma_{N,1} } \bar{\bm t}_n \cdot \bm u \, \mathrm{d}s,
\label{RegVAF0}
\end{equation}
with $\bm u:\Omega\rightarrow\mathbb{R}^d$ and $\alpha:\Omega\rightarrow[0,1]$ standing for the smeared counterparts of the discontinuous displacement and the crack set in (\ref{VAF}). The phase-field variable $\alpha$ takes the value $1$ on $\Gamma_c$, decays smoothly to $0$ in a subset of $\Omega\backslash\Gamma_c$ and then takes the $0$-value in the rest of the domain. With this definition, the limits $\alpha=1$ and $\alpha=0$ represent the fully broken and the intact (undamaged) material phases, respectively, whereas the intermediate range $\alpha\in(0,1)$ mimics the transition zone between them. The function ${\sf g}$ is responsible for the material stiffness degradation. The function ${\sf w}$ defines the decaying profile of $\alpha$, whereas the parameter $0<\ell\ll\mathrm{diam}(\Omega)$ controls the size of the support\footnote{For a continuous real-valued function $f:\Omega\rightarrow\mathbb{R}$, by the support of $f$ we mean all ${\bf x}\in\Omega:f>0$.} of $\alpha$, in other words, the thickness of the transition zone between the two material states.
  
The functions ${\sf g}$ and ${\sf w}$ are the major ingredients of the regularized elastic energy and fracture surface energy functionals, and their specific choice establishes the rigorous link between (\ref{VAF}) and (\ref{RegVAF0}) via the notion of $\Gamma$-convergence, see e.g.\ Braides \cite{Braides1998}, Chambolle \cite{Chambolle2004}, 
also giving a meaning to the induced constant $c_{\sf w}$. Thus, ${\sf g}$ must be a continuous monotonic function that fulfills the properties: ${\sf g}(0)=1$, ${\sf g}(1)=0$, ${\sf g}^\prime(1)=0$ and ${\sf g}^\prime(\alpha)<0$ for $\alpha\in[0,1)$, see e.g.\ Pham et al.\ \cite{Pham2011} for argumentation and discussion. The quadratic polynomial 
\begin{equation}
{\sf g}(\alpha):=(1-\alpha)^2,
\label{g}
\end{equation}
is the simplest choice of the kind. The function ${\sf w}$, also called the local part of the dissipated fracture energy density function \cite{Pham2011}, must be continuous and monotonic such that ${\sf w}(0)=0$, ${\sf w}(1)=1$ and ${\sf w}^\prime(\alpha)\geq0$ for $\alpha\in[0,1]$. The constant $c_{\sf w}:=4\int_0^1\sqrt{{\sf w}(t)}\,\mathrm{d}t$ is a normalization constant in the sense of $\Gamma$-convergence. The two suitable candidates for ${\sf w}$ which are widely adopted read
\begin{equation}
 {\sf w}(\alpha): =\left\{
	\begin{tabular}{l}
	 $\alpha$, \\[0.2cm]
	 $\alpha^2$,
	\end{tabular}
	\right.
	\quad \text{such that} \quad
	c_{\sf w}=\left\{
	\begin{tabular}{ll}
		$\frac{8}{3}$, \\[0.2cm]
		$2$.
	\end{tabular}
	\right.
\label{w}
\end{equation}
It should be noted that formulation (\ref{RegVAF0}) combined with the aforementioned choices for ${\sf g}$ and ${\sf w}$ leads to the so-called $\mathtt{AT}$-1 and $\mathtt{AT}$-2 models, see Table \ref{AT}. $\mathtt{AT}$ stands for {\em Ambrosio-Tortorelli} and the corresponding type of regularization, see \cite{AT1990}. The main difference between the two models is that $\mathtt{AT}$-1 provides the existence of an elastic stage before the onset of fracture, whereas using $\mathtt{AT}$-2 the phase-field starts to evolve as soon as the material is loaded, see e.g.\ \cite{Amor2009,Pham2011} for a more detailed explanation.

\begin{table}[!ht]
\caption{\em Ingredients of formulation (\ref{RegVAF0}).}
\centering
      \begin{tabular}{c|c|c}
      \hline
        ${\sf g}$   &  ${\sf w}$ &  name  \\
	\hline
        $(1-\alpha)^2$ &  \begin{tabular}{c} $\alpha$ \\   $\alpha^2$ \end{tabular}  
				     &  \begin{tabular}{c} $\mathtt{AT}$-1 model \\ $\mathtt{AT}$-2 model \end{tabular} \\
	\hline
      \end{tabular}
\label{AT}
\end{table}
\noindent Other representations for ${\sf g}$ and ${\sf w}$ are available in the literature, see e.g.\ \cite{Borden2012diss,Kuhn2015,Sargado2018,Ambati2015,Alessi2015,Borden2016,Wilson2016,Burke2010b,Burke2013}.

The final note regards the elastic strain energy function $\Psi$. Due to the symmetry of $\Psi$ with respect to the variable $\bm u$, formulation (\ref{RegVAF0}) does not distinguish between fracture behavior in tension and compression. In the numerical simulations, this is manifested by the mesh interpenetration inside of the compressed fractured zones, as reported e.g.\ in Bourdin et al.\ \cite[section 3.3]{Bourdin2000}, Del Piero et al.\ \cite[sections 7]{DelPiero2007} and Lancioni and Royer-Carfagni \cite[Sections 4.1--4.3]{Lancioni2009}. In the discrete crack setting, this would be equivalent to compressive interpenetration of the crack faces. One of the proposed remedies for avoiding such a non-physical behavior implies placing $\Psi$ into the context of non-linear (finite) elasticity, as first presented in \cite{DelPiero2007}. To remain within the framework of linear elasticity, the alternative is to break the symmetry by introducing an additive split of $\Psi$ into the so-called 'tensile' and 'compressive' parts $\Psi^+$ and $\Psi^-$, respectively, and enabling the degradation of $\Psi^+$ only. This option is advocated in \cite{Lancioni2009,Amor2009,Freddi2009,Freddi2010,Miehe2010a,Miehe2010b} and yields the following enhanced representation of (\ref{RegVAF0}):
\begin{equation}
\mathcal{E}(\bm u,\alpha) = \int_\Omega
\left[
{\sf g}(\alpha)\Psi^+(\bm\varepsilon(\bm u))+\Psi^-(\bm\varepsilon(\bm u))
\right]
\mathrm{d}{\bf x}
+ \frac{G_c}{c_{\sf w}} \int_\Omega 
\left(\frac{{\sf w}(\alpha)}{\ell}+\ell|\nabla\alpha|^2\right)
\mathrm{d}{\bf x}
-\int_{\Gamma_{N,1} } \bar{\bm t}_n \cdot \bm u \, \mathrm{d}s.
\label{RegVAF1}
\end{equation}
The two widely adopted splits are based respectively on the so-called volumetric-deviatoric decomposition of the strain tensor $\bm\varepsilon$, as independently presented in Amor et al.\ \cite{Amor2009} and Freddi and Royer-Carfagni \cite{Freddi2009} (based on an extension of the split in \cite{Lancioni2009}), and the spectral decomposition of $\bm\varepsilon$ considered in Miehe et al.\ \cite{Miehe2010a,Miehe2010b}. An idea similar to the latter one is also developed in Freddi and Royer-Carfagni \cite[section 3.4]{Freddi2010}. We refer the interested reader to the aforementioned publications and to \cite{Li2016diss}, where options for constructing $\Psi^\pm$ and the related implications are explained.

Formulation (\ref{RegVAF1}) is what is usually referred to as phase-field model of brittle fracture (at least in the engineering literature). Despite the wide employment of the formulation, the $\Gamma$-convergence result that relates (\ref{RegVAF1}) to the original Francfort-Marigo formulation (\ref{VAF}) is, in general, not available. Some particular results have recently been established in \cite{Chambolle2018}.

In this manuscript, we consider (\ref{RegVAF1}) with both the volumetric-deviatoric and spectral splits (for the sake of compactness, we will denote them as $\mathtt{VD}$- and $\mathtt{S}$-split, respectively), and with the functions ${\sf g}$ and ${\sf w}$ given in Table \ref{AT}.

\subsection{Incremental variational formulation (quasi-static evolution)}
With $\mathcal{E}$ defined by (\ref{RegVAF1}), the state of the system at a given loading step $n\geq 1$ is represented by the solution of
\begin{equation}
\mathrm{arg\,min} \{ \mathcal{E}(\bm u,\alpha): \; {\bm u}\in {\bf V}_{\bar{\bm u}_n}, \alpha\in \mathcal{D}_{\alpha_{n-1}} \},
\label{argmin0}
\end{equation}
where
\begin{equation*}
{\bf V}_{\bar{\bm u}_n}:=\{ \bm u\in{\bf H}^1(\Omega): \; 
\bm u=\bm 0 \; \mathrm{on} \; \Gamma_{D,0}, \; \bm u=\bar{\bm u}_n\; \mathrm{on} \; \Gamma_{D,1} \}
\end{equation*}
is the kinematically admissible displacement space with ${\bf H}^1(\Omega):=[H^1(\Omega)]^d$ and $H^1$ denoting the usual Sobolev space, and
\begin{equation*}
\mathcal{D}_{\alpha_{n-1}}:=\{ \alpha\in H^1(\Omega): \; \alpha\geq\alpha_{n-1} \; \mathrm{in} \; \Omega \}
\end{equation*}
is the admissible space for $\alpha$ with $\alpha_{n-1}$ being known from the previous step. The condition $\alpha\geq\alpha_{n-1}$ in $\Omega$ is used to enforce the {\em irreversibility} of the crack phase-field evolution.

Due to the $\alpha\geq\alpha_{n-1}$ requirement, problem (\ref{argmin0}) is a constrained minimization problem and its necessary optimality condition for computing the solution $(\bm u,\alpha)\in {\bf V}_{\bar{\bm u}_n}\times\mathcal{D}_{\alpha_{n-1}}$ is a variational inequality. Its partitioned form reads as
\begin{equation}
\boxed
{
\left\{
\begin{tabular}{l}
${\mathcal E}_{\bm u}(\bm u,\alpha;\bm v)=0, \quad \forall \; \bm v\in{\bf V}_0$,   \\[0.2cm]
${\mathcal E}_\alpha(\bm u,\alpha;\beta-\alpha)\geq0, \quad \forall \; \beta\in\mathcal{D}_{\alpha_{n-1}}$,
\end{tabular}
\right.
}
\label{Weak0}
\end{equation}
see e.g.\ \cite{Burke2010b,Burke2013,Pham2011,Farrell2017}, where ${\mathcal E}_{\bm u}$ and ${\mathcal E}_\alpha$ are the directional derivatives of the energy functional with respect to $\bm u$ and $\alpha$, respectively,
\begin{equation}
{\mathcal E}_{\bm u}(\bm u,\alpha;\bm v):=\int_\Omega
\left[
{\sf g}(\alpha)\frac{\partial\Psi^+}{\partial\bm\varepsilon}(\bm\varepsilon(\bm u))
+\frac{\partial\Psi^-}{\partial\bm\varepsilon}(\bm\varepsilon(\bm u))
\right]
:\bm\varepsilon(\bm v)\,\mathrm{d}{\bf x}
-\int_{\Gamma_{N,1} } \bar{\bm t}_n \cdot \bm v\,\mathrm{d}s,
\label{Eu}
\end{equation}
\begin{equation}
{\mathcal E}_\alpha(\bm u,\alpha;\beta):=\int_\Omega
\left[
{\sf g}^\prime(\alpha)\Psi^+(\bm\varepsilon(\bm u))\beta
+\frac{G_c}{c_{\sf w}}\left( \frac{1}{\ell}{\sf w}^\prime(\alpha)\beta + 2\ell\nabla{\alpha}\cdot\nabla{\beta} \right)
\right] \,\mathrm{d}{\bf x}.
\label{Ed}
\end{equation}
The displacement test space in (\ref{Weak0}) is defined as ${\bf V}_0:=\{ \bm v\in{\bf H}^1(\Omega): 
\bm v=\bm 0 \; \mathrm{on} \; \Gamma_{D,0}\cup\Gamma_{D,1} \}$.

\subsection{Irreversibility constraint and its treatment}
\label{PIC}
The variational inequality ${\mathcal E}_\alpha\geq 0$ in (\ref{Weak0}) that stems from the irreversibility constraint $\alpha\geq\alpha_{n-1}$ requires special solution algorithms, see e.g.\ \cite{Kinder1980,Glow1981} and \cite[section 5]{Burke2010b}. Below, we outline the available options of handling $\alpha\geq\alpha_{n-1}$ which lead to the equality-based formulations. These include the penalized formulation of our interest. The equivalence between the corresponding formulations and the reference one in (\ref{Weak0}) is highlighted.

\paragraph{{\large\ding{172}} 'Crack-set' irreversibility:} This is a version of irreversibility introduced by Bourdin et al. \cite{Bourdin2000}--\cite{Bourdin2008} and also adopted e.g.\ by Burke et al.\ \cite{Burke2010a}, which relies on the notion of a crack set: if at the current loading step $n$ the (crack) set
\begin{equation}
\mathtt{CR}_{n-1}:=\{{\bf x}\in\overline{\Omega}:\; \alpha_{n-1}({\bf x})\geq \mathtt{CRTOL}\},
\label{CRset}
\end{equation}
where $0\ll\texttt{CRTOL}<1$ is a specified threshold, is non-empty, one explicitly sets $\alpha=1$ for all ${\bf x}\in\texttt{CR}_{n-1}$. The corresponding analogue to (\ref{argmin0}) to be solved at step $n$ is
\begin{equation}
\mathrm{arg\,min} \{ \mathcal{E}(\bm u,d): \; {\bm u}\in {\bf V}_{\bar{\bm u}_n}, \alpha\in\{H^1(\Omega):
\alpha|_{\mathtt{CR}_{n-1}}=1\} \},
\label{argmin3}
\end{equation}
and the resulting weak system of equations for $(\bm u,\alpha)$ reads
\begin{equation}
\boxed
{
\left\{
\begin{tabular}{l}
${\mathcal E}_{\bm u}(\bm u,\alpha;\bm v)=0, \quad \forall \; \bm v\in{\bf V}_0$   \\[0.2cm]
${\mathcal E}_\alpha(\bm u,\alpha;\beta)=0, \quad \forall \; \beta\in\{H^1(\Omega):
\beta|_{\mathtt{CR}_{n-1}}=0\}$,
\end{tabular}
\right.
}
\label{Weak3}
\end{equation}
where ${\mathcal E}_{\bm u}$ and ${\mathcal E}_\alpha$ are given by (\ref{Eu}) and (\ref{Ed}), respectively. As noted in \cite{Amor2009}, the present option can be viewed as an {\em approximate}, or, more precisely, {\em relaxed} version of the requirement $\alpha\geq\alpha_{n-1}$ in $\Omega$ since it only enforces irreversibility of a fully developed crack, whereas phase-field patterns with $\alpha({\bf x})<\mathtt{CRTOL}$ for all ${\bf x}\in\Omega$, which from the mechanical standpoint may be viewed as partially damaged regions, are allowed to heal.

\begin{remark}
Various algorithmic treatments of the 'crack-set' irreversibility can be found in Del Piero et al.\ \cite{DelPiero2007}, and in Lancioni and Royer-Carfagni \cite{Lancioni2009}. Its more sophisticated version is considered in Burke et al.\ \cite{Burke2010b,Burke2013}.
\end{remark}

\begin{remark}
In Artina et al.\ \cite{Artina2014,Artina2015} and Gerasimov and De Lorenzis \cite{Gerasimov2016}, the Dirichlet condition $\alpha|_{\mathtt{CR}_{n-1}}=1$ is realized via penalization by introducing into (\ref{RegVAF1}) the functional
\begin{equation*}
P(\alpha;\tau):=\frac{1}{2\tau}\int_{\mathtt{CR}_{n-1}}(1-\alpha)^2 \, \mathrm{d}{\bf x}, \quad 0<\tau\ll 1,
\end{equation*}
thus yielding the penalized counterpart of (\ref{argmin3}) and (\ref{Weak3}).
\end{remark}

The choice of the threshold value $\mathtt{CRTOL}$ is subtle and may have a strong impact on the computational results, as shown in \cite{Burke2010b,Burke2013}. A similar concern applies to the choice of $\tau$ in the corresponding penalized realization.

\paragraph{{\large\ding{173}} 'History field' irreversibility:} In \cite{Miehe2010b}, Miehe and coworkers proposed the idea of enforcing the irreversibility constraint $\alpha\geq\alpha_{n-1}$ implicitly, via the notion of a {\em history}-field. Their major assumption is that the tensile energy $\Psi^+$ can be viewed as the {\em driving force} of the phase-field evolution and, hence, the maximal $\Psi^+$ accumulated within the loading history, denoted as 
\begin{equation}
\mathcal{H}_n({\bf x}):=
\max_{n\geq 1} \{ \mathcal{H}_{n-1}({\bf x}),\Psi^+(\bm\varepsilon(\bm u)) \},
\label{Hist}
\end{equation}
with $\mathcal{H}_0\equiv0$, must guarantee the fulfillment of $\alpha\geq\alpha_{n-1}$. Technically, one substitutes $\mathcal{H}_n$ to $\Psi^+$ in the original $\mathcal{E}_\alpha$ in (\ref{Ed}) such that
\begin{equation}
\widetilde{\mathcal E}_\alpha(\bm u,\alpha;w):=\int_\Omega
\left[
{\sf g}^\prime(\alpha)\mathcal{H}_n \beta
+\frac{G_c}{c_{\sf w}}\left( \frac{1}{\ell}{\sf w}(\alpha) \beta + 2\ell\nabla{\alpha}\cdot\nabla{\beta} \right)
\right] \,\mathrm{d}{\bf x},
\label{EdMod}
\end{equation}
and forms the system for computing the solution $(\bm u,\alpha)\in {\bf V}_{\bar{\bm u}_n}\times H^1(\Omega)$:
\begin{equation}
\boxed
{
\left\{
\begin{tabular}{l}
${\mathcal E}_{\bm u}(\bm u,\alpha;\bm v)=0, \quad \forall \; \bm v\in{\bf V}_0$   \\[0.2cm]
$\widetilde{\mathcal E}_\alpha(\bm u,\alpha;\beta)=0, \quad \forall \; \beta\in H^1(\Omega)$,
\end{tabular}
\right.
}
\label{Weak4}
\end{equation}
where ${\mathcal E}_{\bm u}$ is given by (\ref{Eu}). System (\ref{Weak4}) is composed of equalities and also uses unconstrained spaces for $\alpha$ and $\beta$. It should be noted, however, that the constructed $\widetilde{\mathcal E}_\alpha$ is no longer of variational nature. Furthermore, the equivalence between (\ref{Weak4}) and the original formulation in (\ref{Weak0}) is not evident and, to the best of our knowledge, neither theoretical nor numerical results that prove it are available. Also, the argumentation from \cite{Miehe2010b} that views $\mathcal{H}_n$ as the driving force appears open to question\footnote{As a first step of this argumentation, one formally ignores the condition $\alpha\geq\alpha_{n-1}$. This enables one to replace the optimality condition ${\mathcal E}_\alpha(\bm u,\alpha;\beta-\alpha)\geq0$ by the simpler, equality-based one, namely, 
\begin{equation}
{\mathcal E}_\alpha(\bm u,\alpha;\beta)=0.
\label{driv0}\tag{A}
\end{equation}
Using this, one arrives at the equation for $\alpha$, which in a strong form reads
\begin{equation*}
-\ell^2\Delta\alpha+\frac{1}{2}{\sf w}^\prime(\alpha)
=-\frac{c_{\sf w}}{2}\frac{\ell}{G_c}{\sf g}^\prime(\alpha)\Psi^+(\bm\varepsilon(\bm u)).
\end{equation*}
With ${\sf g}(\alpha)=(1-\alpha)^2$, ${\sf w}(\alpha)=\alpha^2$ and $c_{\sf w}=2$ considered in \cite{Miehe2010b}, the above turns into
\begin{equation}
-\ell^2\Delta\alpha+\alpha
=\frac{2\ell}{G_c}(1-\alpha)\Psi^+(\bm\varepsilon(\bm u)).
\label{driv1}\tag{B}
\end{equation}
$\Psi^+$ in the right-hand side of (\ref{driv1}) is then interpreted as a source term and, hence, as a driving force of the evolution of $\alpha$. This, in turn, gives rise to the notion of $\mathcal{H}_n$ in (\ref{Hist}) and to the claim that the introduction of $\mathcal{H}_n$ into (\ref{driv1}) (equivalently, into (\ref{driv0})) will enforce $\alpha\geq\alpha_{n-1}$. It is easy to see, however, that re-arranging (\ref{driv1}) yields
\begin{equation*}
-\ell^2\Delta\alpha+\left(1+\frac{2\ell}{G_c}\Psi^+(\bm\varepsilon(\bm u))\right)\alpha
=\frac{2\ell}{G_c}\Psi^+(\bm\varepsilon(\bm u)),
\end{equation*}
where $\Psi^+$ enters also the left-hand-side of the equation as a coefficient by the linear term. Thus, the role of $\Psi^+$ as a driving force is not evident.
}. Finally, it can be anticipated -- see sections \ref{IntroMR} and \ref{Basics} -- that in the case ${\sf w}(\alpha)=\alpha$, formulation (\ref{Weak4}) will not work without the additionally enforced condition $\alpha\geq0$.

Interestingly, in the earlier paper of Miehe et al.\ \cite{Miehe2010a}, the authors considered the penalization option similar to representation (\ref{OurPen}) below, yet already in \cite{Miehe2010b} they switched to the notion of $\mathcal{H}_n$ and have been using it in all their following publications.

\paragraph{{\large\ding{174}} Penalization:}
The most straightforward way of addressing $\alpha\geq\alpha_{n-1}$, which results in the {\em equality} to be solved instead of the inequality, is via penalization: one should add the penalty term
\begin{equation}
P(\alpha;\gamma):=\frac{\gamma}{2}\int_\Omega \langle\alpha-\alpha_{n-1} \rangle_{-}^2 \; \mathrm{d}{\bf x}, 
\quad \gamma\gg 1,
\label{OurPen}
\end{equation}
to the energy functional $\mathcal{E}$ in (\ref{RegVAF1})\footnote{A general form of the integrand in (\ref{OurPen}) is $\langle\alpha-\alpha_{n-1} \rangle_{-}^p$ with $p\geq1$ and $P$ represents a regularization of the indicator function 
\begin{equation*}
I_{\mathcal{D}_{\alpha_{n-1}}}(\alpha):=
\left\{
\begin{tabular}{cl}
$0$, & in $\mathcal{D}_{\alpha_{n-1}}$, \\ [0.1cm]
$+\infty$, & otherwise,
\end{tabular}
\right.
\label{indicator}\tag{C}
\end{equation*}
to be originally added to $\mathcal E$ in (\ref{RegVAF1}).}. 
In (\ref{OurPen}), $\langle y \rangle_{-}:=\min(0,y)$. The corresponding variational problem reads
\begin{equation}
\mathrm{arg\,min} \{ \mathcal{E}(\bm u,\alpha)+P(\alpha;\gamma): \; {\bm u}\in {\bf V}_{\bar{\bm u}_n}, \alpha\in H^1(\Omega) \}.
\label{argmin1}
\end{equation}
and the resulting weak system for $(\bm u,\alpha)$ is as follows
\begin{equation}
\boxed
{
\left\{
\begin{tabular}{l}
${\mathcal E}_{\bm u}(\bm u,\alpha;\bm v)=0, \quad \forall \; \bm v\in{\bf V}_0$,   \\[0.2cm]
$\displaystyle {\mathcal E}_\alpha(\bm u,\alpha;\beta)+\gamma\int_\Omega\langle \alpha-\alpha_{n-1} \rangle_{-}\beta=0, \quad \forall \; \beta\in H^1(\Omega)$,
\end{tabular}
\right.
}
\label{Weak1}
\end{equation}
with ${\mathcal E}_{\bm u}$ and ${\mathcal E}_\alpha$ given again by (\ref{Eu}) and (\ref{Ed}), respectively. Note that the admissible space for $\alpha$ and the test space for $\beta$ are no longer constrained. The obtained penalized un-constrained problem (\ref{Weak1}) approximates the original constrained problem (\ref{Weak0}) and thus is {\em equivalent} to it in the limit of $\gamma\rightarrow\infty$. Thus the appropriate choice of $\gamma$ is always viewed as a critical point of the technique also within the numerical experiments. A too small penalty parameter will lead to inaccurate (i.e.\ insufficient) enforcement of the constraint, a too large one will result in ill-conditioning.

In section \ref{MainRes}, we devise an analytical procedure for the 'reasonable' choice of a {\em lower bound} for $\gamma$ that guarantees a sufficiently accurate enforcement of the crack phase-field irreversibility constraint. Furthermore, in our numerical studies the aforementioned findings will be validated. We will also compare the solutions of problems (\ref{Weak4}) and (\ref{Weak1}) for a benchmark case.

\begin{remark}
Another constraint enforcement technique, namely, the augmented Lagrangian approach is adopted by Wheeler et al.\ \cite{Wheeler2014}, and Wick \cite{Wick2016,Wick2017}. It implies adding to (\ref{RegVAF1}) the functional
\begin{equation}
P(\alpha;\Xi,\gamma):=\frac{1}{2\gamma}\int_\Omega
\langle\Xi + \gamma\alpha \rangle_{+}^2 \; \mathrm{d}{\bf x}
+\frac{1}{2\gamma}\int_\Omega
\langle\Xi + \gamma(\alpha-\alpha_{n-1}) \rangle_{-}^2 \; \mathrm{d}{\bf x}
-\frac{1}{2\gamma}\int_\Omega\Xi^2 \; \mathrm{d}{\bf x},
\label{WickPen}
\end{equation}
where $\Xi\in L^2(\Omega)$ is an unknown field and, again, $\gamma\gg 1$ is a user-prescribed penalty constant. The above $P$ is the so-called Moreau-Yosida approximation of the indicator function $I_{\mathcal{D}_{\alpha_{n-1}}}$ in (\ref{indicator}), see \cite{Wheeler2014} for details. It is claimed that using (\ref{WickPen}), the possible stability issues occurring in case of (\ref{OurPen}) are avoided. However, the appropriate choice of $\gamma$ remains a similar issue as in (\ref{OurPen}), and, furthermore, with the extra variable $\Xi$ to be solved for, the computational effort is significantly increased. 
\end{remark}

\subsection{Staggered solution scheme}
\label{StagSol}
The staggered solution algorithm for the system in (\ref{Weak3}), (\ref{Weak4}) and (\ref{Weak1}) implies alternately fixing $\bm u$ and $\alpha$, and solving the corresponding equations until convergence. The algorithm is sketched in Table \ref{TableStag} with the step 2 specifically adjusted to the evolution equation in (\ref{Weak1}). Adaptation of this step w.r.t.\ the corresponding equation in (\ref{Weak3}), (\ref{Weak4}) is straightforward\footnote{In particular, in the case of (\ref{Weak4}), the quantity $\mathcal{H}_n$ entering $\widetilde{\mathcal{E}}_\alpha$ must be re-defined: at every $k\geq1$ we set $\mathcal{H}_n:=\max\{\mathcal{H}_{n-1},\Psi^+(\bm\varepsilon(\bm u^{k-1})) \}$.}. 

\begin{table}[!ht]\footnotesize
\caption{\em Staggered iterative solution process for (\ref{Weak1}) at a fixed loading step $n\geq 1$.}
\label{TableStag}
{
\begin{tabular}{l}
\hline \\
{\bf Input:} loading data $(\bar{\bm u}_n,\bar{\bm t}_n)$ on $\Gamma_{D,1},\Gamma_{N,1}$, and \\ 

{\color{white}\bf Input:} solution $(\bm u_{n-1},\alpha_{n-1})$ from step $n-1$. \\ [0.1cm]

\quad\quad Initialization, $k=0$:\\
\quad\quad\quad 1.\; set $(\bm u^0,\alpha^0):=(\bm u_{n-1},\alpha_{n-1})$.\\ [0.1cm]

\quad\quad Staggered iteration $k\geq1$:\\

\quad\quad\quad 2.\; given $\bm u^{k-1}$, solve $\displaystyle {\mathcal E}_\alpha(\bm u^{k-1},\alpha;\beta)+
\gamma\int_\Omega\langle \alpha-\alpha_{n-1} \rangle_{-}\beta=0$ for $\alpha$, set $\alpha=:\alpha^k$,\\

\quad\quad\quad 3.\; given $\alpha^k$, solve ${\mathcal E}_{\bm u}(\bm u,\alpha^k;\bm v)=0$ for $\bm u$, set $\bm u=:\bm u^k$,\\

\quad\quad\quad 4.\; for the obtained pair $(\bm u^k,\alpha^k)$, check \\

\quad\quad\quad\quad\quad\quad\quad
$\mathrm{Res}_\mathrm{Stag}^k:=|\displaystyle {\mathcal E}_\alpha(\bm u^k,\alpha^k;\beta)+
\gamma\int_\Omega\langle \alpha^k-\alpha_{n-1} \rangle_{-}\beta|
\leq\texttt{TOL}_\mathrm{Stag}, \; \forall \; \beta\in H^1(\Omega)$, \\

\quad\quad\quad 5.\; if fulfilled, set $(\bm u^k,\alpha^k)=:(\bm u_n,\alpha_n)$ and stop; \\

\quad\quad\quad {\color{white}5.}\; else $k+1\rightarrow k$. \\ [0.1cm]

{\bf Output:} solution $(\bm u_n,\alpha_n)$. \\ \\
\hline
\end{tabular}
}
\end{table}

Both equations in Table \ref{TableStag} are non-linear: for the first one this is due to the non-linearity of ${\sf g}(\alpha)\frac{\partial\Psi^+}{\partial\bm\varepsilon}(\bm\varepsilon(\bm u))$ $+\frac{\partial\Psi^-}{\partial\bm\varepsilon}(\bm\varepsilon(\bm u))=:\bm\sigma(\bm u,\alpha)$, whereas for the second one it is due to the Macaulay bracket term. Therefore, the Newton-Raphson procedure is used to iteratively compute $\bm u^k$ and $\alpha^k$ with $\bm u^{k-1}$ and $\alpha^{k-1}$ being taken as the corresponding initial guesses, and $\mathtt{TOL}_\mathrm{NR}$ as the tolerance. Owing to the nested nature of the Newton-Raphson loops, we impose that $\mathtt{TOL}_\mathrm{NR}<\mathtt{TOL}_\mathrm{Stag}$.

The definition of $\mathrm{Res}_\mathrm{Stag}^k$ in Table \ref{TableStag} used for stopping the staggered process is not unique. As an example, the quantity $||\alpha^k-\alpha^{k-1}||_\infty$ is considered as the residual e.g.\ in \cite{Amor2009,Pham2011,Mesgarnejad2015,Farrell2017} when solving (\ref{Weak0}) and in \cite{Bourdin2000,Bourdin2007a,Bourdin2007b,Bourdin2008} while solving (\ref{Weak3}). Another option for $\mathrm{Res}_\mathrm{Stag}^k$ is a relative change in the normalized energy $\mathcal{E}$. It has been taken as a residual in \cite{Ambati2015review} while solving (\ref{Weak4}), and as an auxiliary convergence-tracing quantity in \cite{Gerasimov2016} while considering (\ref{Weak3}). 

Finally, as already mentioned in the introductory part, the problem size of the system in (\ref{Weak0}), (\ref{Weak3}), (\ref{Weak1}) and (\ref{Weak4}) after finite element discretization is typically very large, since both the phase-field and the deformation localize in bands of width of order $\ell$. Solving the system in a staggered way in this case is computationally very demanding, see e.g.\ \cite{Mesgarnejad2015,Ambati2015review,Gerasimov2016,Farrell2017} for detailed studies.

\section{Penalized formulations}
\label{MainRes}
\subsection{Outline of the procedure and main results}
\label{IntroMR}
In this section, a procedure for the optimal choice of the penalty constant $\gamma\gg 1$ in formulation (\ref{argmin1}) is devised. For better clarity, we briefly outline the path followed and the main results.

\begin{figure}[!ht]
\begin{center}
\includegraphics[width=1.0\textwidth]{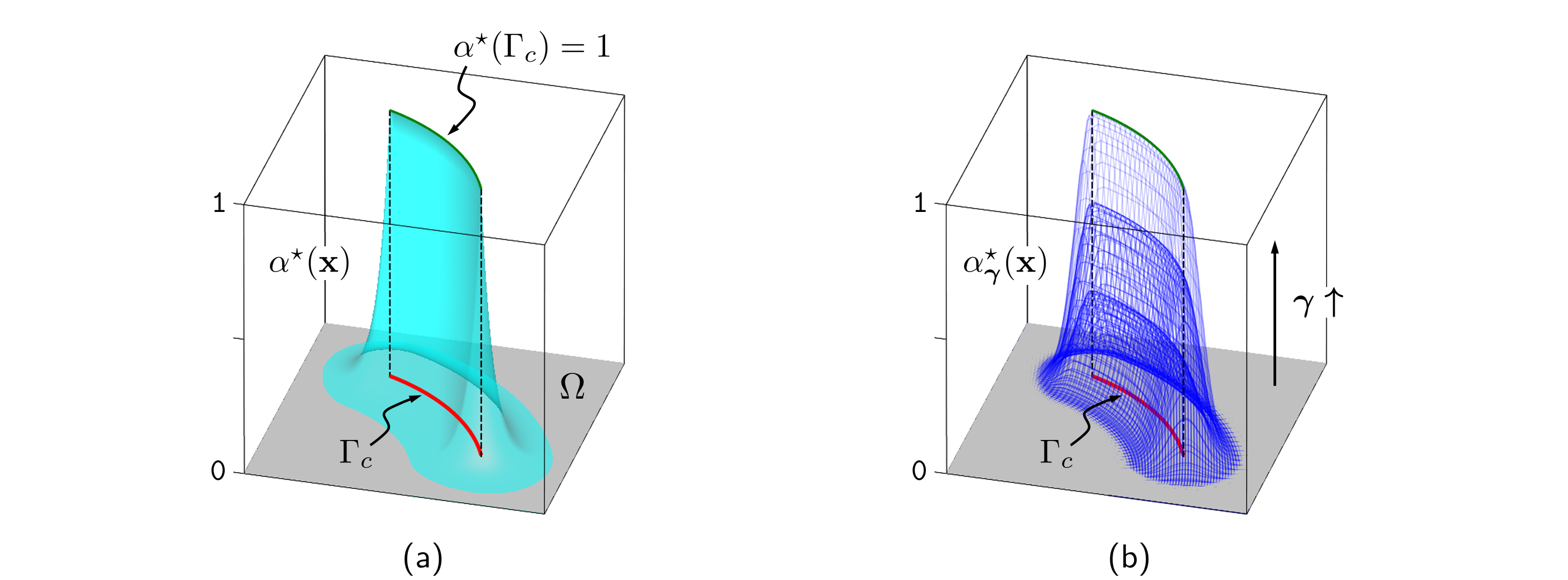}
\end{center}
\caption{Sketch of the optimal phase-field profile that corresponds to a pre-existing or fully-developed crack $\Gamma_c$ in two dimensions: (a) $\alpha^\star$, minimizer of $E_S$ in (\ref{Es}); (b) $\alpha^\star_\gamma$, minimizer of $\widetilde{E}_S$ in (\ref{Es_gam}), as a function of the penalty parameter $\gamma>0$.}
\label{Fig2}
\end{figure}

Our starting point is the notion of the {\em optimal phase-field profile}, as well as the related $\Gamma$-convergence result available for the regularized {\em fracture surface energy} functional $E_S$ in the original representations (\ref{RegVAF0}) and (\ref{RegVAF1})
\begin{equation}
E_S(\alpha):=
\frac{G_c}{c_{\sf w}} \int_\Omega 
\left(\frac{{\sf w}(\alpha)}{\ell}+\ell|\nabla\alpha|^2\right)\mathrm{d}{\bf x},
\label{Es}
\end{equation}
The $\Gamma$-convergence framework states the following: if $\Gamma_c\subset\Omega$ is a pre-existing, or a fully-developed crack, a phase-field profile $\alpha^\star:\Omega\rightarrow[0,1]$ that corresponds to $\Gamma_c$ -- in the literature it is usually termed {\em optimal} -- is a solution of the minimization problem for $E_S$ on a suitable admissible set. More precisely, it is 
\begin{equation}
\alpha^\star:=\mathrm{arg\,min}
\left\{ E_S(\alpha), \;
\alpha\in H^1(\Omega): \, \alpha(\Gamma_c)=1, \,
\alpha\geq0 \right\},
\label{aStar}
\end{equation}
see Figure \ref{Fig2}(a) for an illustrative sketch in a two-dimensional setting. The given $\alpha^\star$ is such that
\begin{equation}
G_c^{-1}E_S(\alpha^\star)=|\Gamma_c|, \;\; \mathrm{if} \;\; {\sf w}(\alpha)=\alpha,
\label{GconvLin}
\end{equation}
and
\begin{equation}
\lim_{\ell\rightarrow 0} G_c^{-1}E_S(\alpha^\star)=|\Gamma_c|, \;\; \mathrm{if} \;\; {\sf w}(\alpha)=\alpha^2.
\label{GconvSq}
\end{equation}
Note that if $0<\ell\ll\mathrm{diam}(\Omega)$ is finite, one expects in the latter case that $G_c^{-1}E_S(\alpha^\star)\approx|\Gamma_c|$. In the one-dimensional setting, $\alpha^\star$ in (\ref{aStar}) and the related $\Gamma$-convergence property\footnote{In higher dimensions $|\Gamma_c|=\int_{\Gamma_c} \mathrm{d}s=\int_\Omega \delta({\bf x}-\Gamma_c)\, \mathrm{d}{\bf x}$, with $\delta$ as the Dirac delta. In one dimension, assuming $\Omega:=(-L,L)$ such that $\Gamma_c$ is represented by the point $x=0$, it is $|\Gamma_c|=\int_{-L}^L \delta(x)\, \mathrm{d}x=1$.} can be constructed explicitly, as detailed in section \ref{Basics}.

The above enables us to similarly address in section \ref{Penal1} the proposed penalized formulation (\ref{argmin1}). For this, we assume first that the penalty term $P$ given by (\ref{OurPen}) must be included into the fracture surface energy, thus yielding the following modification of (\ref{Es}):
\begin{equation}
\widetilde{E}_S(\alpha):=
\frac{G_c}{c_{\sf w}} \int_\Omega 
\left(\frac{{\sf w}(\alpha)}{\ell}+\ell|\nabla\alpha|^2\right)\mathrm{d}{\bf x}
+\frac{\gamma}{2}\int_\Omega\langle \alpha-\alpha_{n-1} \rangle_{-}^2 \; \mathrm{d}{\bf x}.
\label{Es_gam}
\end{equation}
In this case, we look for the optimal phase-field profile reading
\begin{equation}
\alpha^\star_\gamma:=\mathrm{arg\,min} \{ \widetilde{E}_S(\alpha), \;
\alpha\in H^1(\Omega) \}.
\label{aStar_gam}
\end{equation}
The subscript indicates that $\alpha^\star_\gamma$ is a parametric function of $\gamma>0$, see the sketch in Figure \ref{Fig2}(b). Assuming $\alpha_{n-1}$ represents a fully developed crack, the Dirichlet boundary condition $\alpha(\Gamma_c)=1$ on $\Gamma_c$ and $\alpha\geq0$ which enter (\ref{aStar}) are naturally omitted in the present admissible space.

Our main result to be derived in section \ref{Penal1} is that for $\alpha^\star_\gamma$ defined by (\ref{aStar_gam}), the $\Gamma$-convergence result holds. Regardless whether ${\sf w}$ is linear or quadratic, it takes the following form:
\begin{equation}
\lim_{\gamma\rightarrow \infty} 
\left( \lim_{\ell\rightarrow 0} G_c^{-1}\widetilde{E}_S(\alpha^\star_\gamma) \right)
=|\Gamma_c|.
\label{Gconv_gam}
\end{equation}
The explicit construction of $\alpha^\star_\gamma$ and the proof of (\ref{Gconv_gam}) is straightforward to give in a one-dimensional setting. Introducing the notation
\begin{equation}
\lim_{\ell\rightarrow 0} G_c^{-1}\widetilde{E}_S(\alpha^\star_\gamma)=:F(\gamma),
\label{F_gam}
\end{equation}
and using (\ref{Gconv_gam}), the range of $\gamma>0$ which we term optimal will finally be derived. It is a solution to
\begin{equation}
(1-\mathtt{TOL}_\mathrm{ir})|\Gamma_c| \leq F(\gamma)<|\Gamma_c|,
\label{gam_opt}
\end{equation}
where $0<\mathtt{TOL}_\mathrm{ir}\ll 1$ is a user-prescribed irreversibility tolerance. In particular, the solution of $F(\gamma)\overset{!}{=}(1-\mathtt{TOL}_\mathrm{ir})|\Gamma_c|$ yields a lower bound for $\gamma$ that guarantees the prescribed accuracy with respect to the reference value $|\Gamma_c|$. This clarifies the meaning of the optimal range we have heuristically introduced above. Interestingly, by assessing the first derivative of $F$ in (\ref{F_gam}), it can also be seen in section \ref{Penal1} that no improvement in terms of the $\Gamma$-convergence recovery is achieved when taking $\mathtt{TOL}_\mathrm{ir}<0.01$. We then can argue that $\mathtt{TOL}_\mathrm{ir}:=0.01$ can be treated as the sufficient practical threshold value. That is, the choice of the threshold value is not arbitrary.

\begin{remark}
The constraint $\alpha\geq0$ present in the admissible space in (\ref{aStar}) must only be introduced for the case of $E_S$ in (\ref{Es}) that uses $\sf{w}(\alpha)=\alpha$. Lacking this, as shown in appendix \ref{Ap0}, the corresponding $\alpha^\star$ becomes negative in a sub-domain of $\Omega$. In contrast, the minimizer of $E_S$ using $\sf{w}(\alpha)=\alpha^2$ is automatically in the range of $[0,1]$, see section \ref{Basics}.
\end{remark}

Following this remark, we finally consider in section \ref{Penal2} a penalized way of treating the constraint $\alpha\geq0$ in the corresponding case. For that, the following modification of $E_S$ is considered (with ${\sf w}(\alpha)=\alpha$ and $c_{\sf w}=\frac{8}{3}$):    
\begin{equation}
\widehat{E}_S(\alpha):=
\frac{3}{8}G_c \int_\Omega 
\left(\frac{\alpha}{\ell}+\ell|\nabla\alpha|^2\right)\mathrm{d}{\bf x}
+\frac{\rho}{2}\int_\Omega\langle \alpha \rangle_{-}^2 \; \mathrm{d}{\bf x},
\label{Es_rho}
\end{equation}
with $\rho\gg 1$, such that the minimizer (optimal profile) to be sought is
\begin{equation}
\alpha^\star_\rho:=\mathrm{arg\,min}
\{ \widehat{E}_S(\alpha), \;
\alpha\in H^1(\Omega): \, \alpha(\Gamma_c)=1 \}.
\label{aStar_rho}
\end{equation}
The procedure for deriving the optimal range of the penalty parameter $\rho$ is then devised. We employ a very similar logic as proposed in section \ref{Penal1}.

As follows, we term $E_S$ (as well as $\widetilde{E}_S$) with ${\sf w}(\alpha)=\alpha$ and ${\sf w}(\alpha)=\alpha^2$ the {\em linear} and {\em quadratic} models, respectively. With this nomenclature, $\widehat{E}_S$ will be referred to as the linear model. Also, we will refer to penalization in the minimization problems for $\widetilde{E}_S$ in (\ref{Es_gam}) and $\widehat{E}_S$ in (\ref{Es_rho}) as the $\gamma$- and $\rho$-penalization, respectively.

Our final note regards the conceptual difference between the two penalized cases (\ref{Es_gam}) and (\ref{Es_rho}). The $\gamma$-penalization aims at handling the {\em irreversibility} during the crack phase-field evolution (regardless of the model type). The $\rho$-penalization is intended for the {\em recovery} of the bounded profile corresponding to a pre-existing crack in case of the linear model.

\subsection{Optimal phase-field profile for $E_S$}
\label{Basics}
In this section, we provide the necessary basic machinery for constructing the optimal profile $\alpha^\star$ for $E_S$ in (\ref{Es}) reduced to the one-dimensional setting. The linear model case is not trivial, since one deals with variational and differential inequalities. We were not able to find in the literature any related results for this case and, therefore, they are elaborated in detail in appendix \ref{Ap0}. The corresponding result for the quadratic model is well-documented in the literature and we only briefly recall it. The $\Gamma$-convergence property for $\alpha^\star$ in either model case is verified as well.

Without loss of generality, we assume $\Omega:=(-L,L)$ such that $\Gamma_c$ is represented by the point $x=0$. Due to symmetry, we write
\begin{equation}
E_S(\alpha):=
2\frac{G_c}{c_{\sf w}} \int_0^L 
\left(\frac{{\sf w}(\alpha)}{\ell}+\ell\left( \alpha^\prime \right)^2\right)\mathrm{d}x.
\label{Es1d}
\end{equation}
The strong formulation for obtaining the minimizer $\alpha^\star$ is as follows
\begin{equation}
\left\{
\begin{tabular}{rl}
$\displaystyle -\ell^2\alpha^{\prime\prime}+\frac{1}{2}\geq0$ & in $(0,L)$,   \\ 
$\alpha=1$ & at $x=0$, \\
$\alpha^\prime\geq0$ & at $x=L$.
\end{tabular}
\right.
\label{bvp1d_Lin}
\end{equation}
when ${\sf w}(\alpha)=\alpha$, and
\begin{equation}
\left\{
\begin{tabular}{rl}
$\displaystyle -\ell^2\alpha^{\prime\prime}+\alpha =0$ & in $(0,L)$,   \\ 
$\alpha=1$ & at $x=0$, \\
$\alpha^\prime=0$ & at $x=L$.
\end{tabular}
\right.
\label{bvp1d_Sq}
\end{equation}
when ${\sf w}(\alpha)=\alpha^2$.

\begin{figure}[!ht]
\begin{center}
\includegraphics[width=1.0\textwidth]{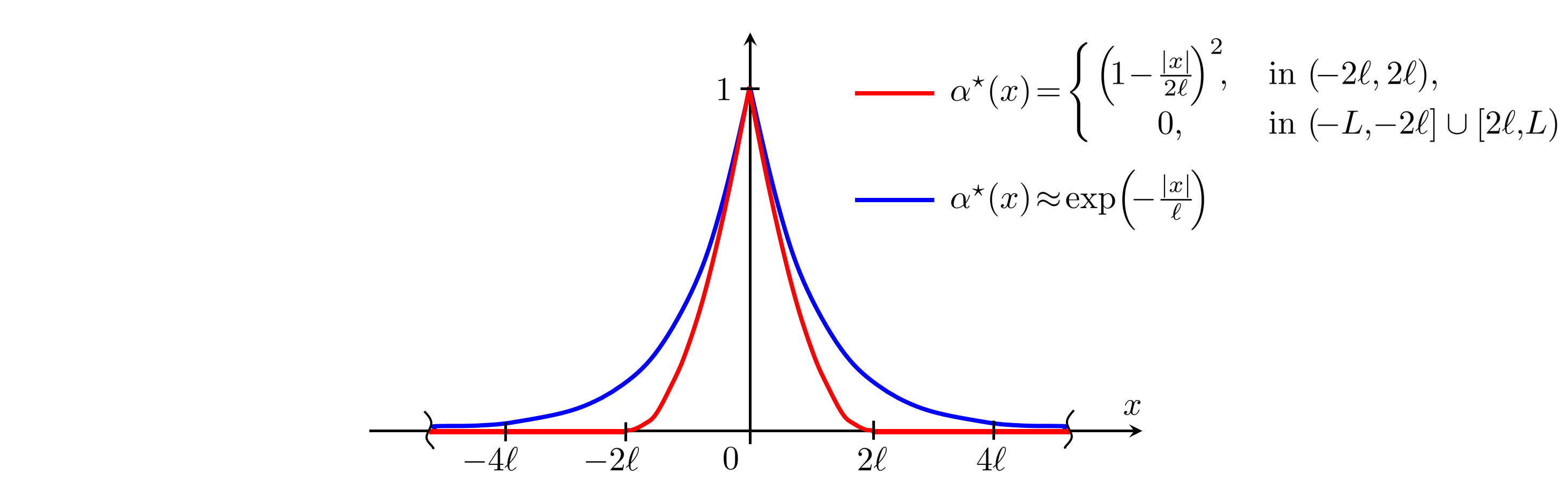}
\end{center}
\caption{Sketch of the optimal profile $\alpha^\star$ in $(-L,L)$ for the linear and quadratic models.}
\label{Fig3}
\end{figure}

As shown in appendix \ref{Ap0}, the admissible set of solutions to (\ref{bvp1d_Lin}) is given by  
\begin{equation*}
0\leq \alpha(x) \leq \left(1- \frac{x}{2\ell}\right)^2 \quad \mathrm{in} \; (0,L).
\end{equation*}
The unique representative $\alpha^\star$ which complies with the crack phase-field formalism in the corresponding interval is a piecewise function which is continuously-differentiable at $x=2\ell$ reading 
\begin{equation}
\alpha^\star(x)=
\left\{
\begin{tabular}{cl}
$\displaystyle
\left(1- \frac{x}{2\ell}\right)^2$, & in $(0,2\ell)$,   \\[0.2cm]
$0$, & in $[2\ell,L)$.
\end{tabular}
\right.
\label{aStar1d_Lin}
\end{equation}
It is straightforward that $G_c^{-1}E_S(\alpha^\star)=1$ regardless of $\ell$, or $\frac{L}{\ell}$. The symmetry argument finally yields the representation which is valid in $(-L,L)$, see Figure \ref{Fig3}. It is bounded within $[0,1]$ and has finite support. 

Solution to problem (\ref{bvp1d_Sq}) is as follows:
\begin{equation}
\alpha^\star(x)=\frac{1}{\exp\left(2\frac{L}{\ell}\right)+1} \exp\left(\frac{x}{\ell}\right)
			  +\frac{\exp\left(2\frac{L}{\ell}\right)}{\exp\left(2\frac{L}{\ell}\right)+1} \exp\left(-\frac{x}{\ell}\right).
\label{aStar1d_Sq}
\end{equation}
It is such that
\begin{equation*}
G_c^{-1}E_S(\alpha^\star)=\frac{\exp\left(2\frac{L}{\ell}\right)-1}{\exp\left(2\frac{L}{\ell}\right)+1},
\end{equation*}
and with the assumption $\ell\ll L$, it holds that $\exp\left(2\frac{L}{\ell}\right)\gg1$, thus yielding $G_c^{-1}E_S\approx1$. The resulting simplified form of the obtained $\alpha^\star$, already known in the literature, follows:
\begin{equation*}
\alpha^\star(x)\approx\exp\left(-\frac{x}{\ell}\right).
\end{equation*}
The symmetry argument ultimately yields the representation which is valid in $(-L,L)$, see Figure \ref{Fig3}. The obtained $\alpha^\star$ is bounded by the range $[0,1]$ and has infinite support.

The interesting implication of the above one-dimensional results and their generalization to higher dimensions is as follows. The minimizer of the linear model has a {\em finite} support within $\Omega$, and therefore the $\Gamma$-convergence property (\ref{GconvLin}) is independent on $\ell$. In contrast, the support of the minimizer of the quadratic model is the entire domain $\Omega$, thus resulting in (\ref{GconvSq}). In other words, the linear model is capable of reproducing the fracture energy $E_S$ in (\ref{VAF}) exactly, whereas the quality of the corresponding approximation in the case of the quadratic model strongly depends on $\ell$.

\subsection{$\gamma$-penalization: irreversibility for the linear and quadratic models}
\label{Penal1}
The one-dimensional counterpart of $\widetilde{E}_S$ in (\ref{Es_gam}) reads
\begin{equation}
\widetilde{E}_S(\alpha):=
2\frac{G_c}{c_{\sf w}} \int_0^L
\left(\frac{{\sf w}(\alpha)}{\ell}+\ell \left( \alpha^\prime \right)^2\right)\mathrm{d}x
+\gamma\int_0^L\langle \alpha-\alpha_{n-1} \rangle_{-}^2 \; \mathrm{d}x.
\label{Es1d_gam}
\end{equation}
(Note again that owing to symmetry we consider the half-length representation.) As mentioned in section \ref{IntroMR}, we will construct the minimizer $\alpha^\star_\gamma$ of (\ref{Es1d_gam}), establish the related $\Gamma$-convergence property and, based on that, derive the optimal range for the penalty constant $\gamma\gg1$ in (\ref{Es1d_gam}). 

The weak formulation is to find $\alpha\in H^1(0,L)$ such that
\begin{equation}
\widetilde{E}_S^\prime(\alpha,\beta)=2\frac{G_c}{c_{\sf w}}
\int_0^L \left(
\frac{{\sf w}^\prime(\alpha)}{\ell}\beta+2\ell \alpha^\prime\beta^\prime
\right) \mathrm{d}x
+2\gamma \int_0^L \langle \alpha-\alpha_{n-1} \rangle_{-}\beta\; \mathrm{d}x
\overset{!}{=}0, \quad \forall \beta\in H^1(0,L).
\label{weak_gam}
\end{equation}
For deriving (\ref{weak_gam}), we recall the definition $\langle y \rangle_{-}:=\min(0,y)$ and also use that $\langle y \rangle_{-}^\prime=:{\sf H}^{-}(y)$, with ${\sf H}^{-}$ as a step-function defined as follows:
\begin{equation}
{\sf H}^{-}(y):=\left\{
\begin{tabular}{ll}
$1$ & $y\leq0$,   \\ 
$0$ & $y>0$.
\end{tabular}
\right.
\label{Hmin}
\end{equation}
The resulting boundary value problem for $\alpha^\star_\gamma$ reads
\begin{equation}
\left\{
\begin{tabular}{rl}
$\displaystyle -\ell^2\alpha^{\prime\prime}+\frac{1}{2}{\sf w}^\prime(\alpha) 
+\frac{c_{\sf w}}{2}\frac{\gamma\ell}{G_c} \langle \alpha-\alpha_{n-1} \rangle_{-}=0$ & in $(0,L)$,   \\ 
$\alpha^\prime=0$ & at $x\in\{0,L\}$.
\end{tabular}
\right.
\label{bvp_gam}
\end{equation}

For further analysis, it proves convenient to introduce the following definition:
\begin{equation}
\boxed
{
\frac{c_{\sf w}}{2}\frac{\gamma\ell}{G_c}=:s.
}
\label{s_variable}
\end{equation}
and we call the dimensionless quantity $s>0$ penalty parameter as well. Furthermore, we assume that $\alpha_{n-1}$ in (\ref{bvp_gam}) -- the phase-field profile known from the previous loading step -- already represents the fully-developed crack. That is, we set
\begin{equation*}
 \boxed
{
\alpha_{n-1}:=
\left\{
\begin{tabular}{ll}
$\alpha^\star$ in (\ref{aStar1d_Lin}), & when ${\sf w}(\alpha)=\alpha$,   \\  [0.2cm]
$\alpha^\star$ in (\ref{aStar1d_Sq}),  & when ${\sf w}(\alpha)=\alpha^2$.
\end{tabular}
\right.
}
\end{equation*}
Using the above, (\ref{bvp_gam}) is compactly rewritten as follows:
\begin{equation}
\left\{
\begin{tabular}{rl}
$\displaystyle -\ell^2\alpha^{\prime\prime}+\frac{1}{2}{\sf w}^\prime(\alpha) 
+s \langle \alpha-\alpha^\star \rangle_{-}=0$ & in $(0,L)$,   \\ 
$\alpha^\prime=0$ & at $x\in\{0,L\}$.
\end{tabular}
\right.
\label{bvp_s}
\end{equation}
In either case of linear and quadratic model, problem (\ref{bvp_s}) is a semi-linear problem and requires an iterative solution procedure. The corresponding derivations are detailed in appendix \ref{Ap1} and below we only summarize them.

\subsubsection{Linear model, solution to (\ref{bvp_s})}
\label{LinMod_gam}
The exact solution to problem (\ref{bvp_s}) can be constructed using the Newton-Raphson procedure in only two iterations $i\geq0$, see appendix \ref{Ap1_LinMod}. It reads:
\begin{equation}
\alpha^\star_\gamma(x):=\alpha^\star(x)+\Delta\alpha_0(x)+\Delta\alpha_1(x),
\label{aStarGamm_Lin}
\end{equation}
where $\Delta\alpha_0$ and $\Delta\alpha_1$ are explicitly given in appendix \ref{Ap1_LinMod} by equations (\ref{Da0_Lin}) and (\ref{Da1_Lin}), respectively.

\begin{figure}[!ht]
\begin{center}
\includegraphics[width=1.0\textwidth]{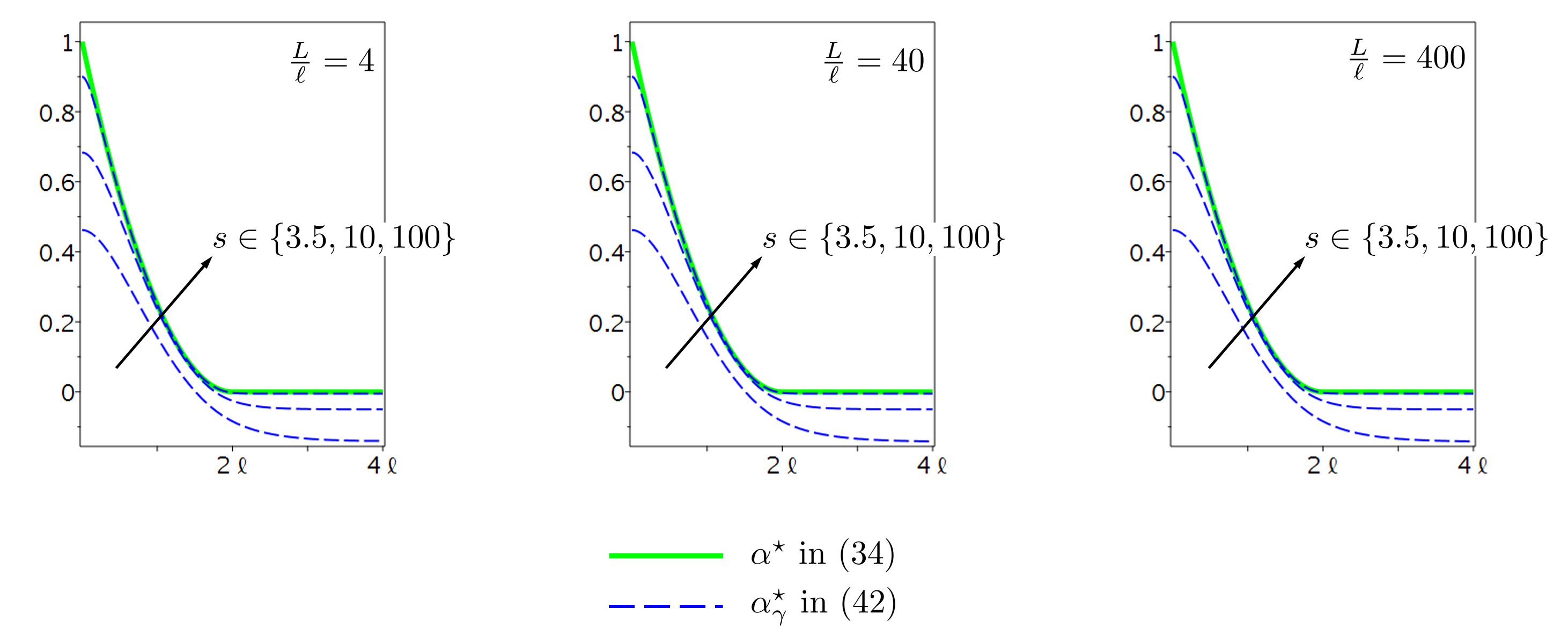}
\end{center}
\caption{The functions $\alpha^\star$ in (\ref{aStar1d_Lin}) and $\alpha^\star_\gamma$ in (\ref{aStarGamm_Lin}) plotted on the interval $x\in[0,4\ell]$ for $\frac{L}{\ell}\in\{4,40,400\}$; for fixed $\frac{L}{\ell}$, $\alpha^\star_\gamma$ is depicted for different values of the penalty parameter $s>0$.}
\label{Fig4}
\end{figure}

In Figure \ref{Fig4}, we plot and compare the profile $\alpha^\star_\gamma$ in (\ref{aStarGamm_Lin}) with the reference representation $\alpha^\star$ in (\ref{aStar1d_Lin}) for three different magnitudes of $\frac{L}{\ell}>2$. Also, for a fixed $\frac{L}{\ell}$, $\alpha^\star_\gamma$ is depicted for different values of the penalty parameter. Note that for the considered cases $\frac{L}{\ell}\in\{4,40,400\}$, both $\alpha^\star$ and $\alpha^\star_\gamma$ are defined, respectively, on the intervals $x\in[0,4\ell]$, $[0,40\ell]$ and $[0,400\ell]$, yet for comparison purposes all plots are restricted to the interval $[0,4\ell]$. First conclusion to be drawn from the figure is as follows: for fixed $\frac{L}{\ell}$ and varying $s$, $\alpha^\star_\gamma$ demonstrates the desired trend, namely, the point-wise convergence of $\alpha^\star_\gamma$ to $\alpha^\star$ when $s$ increases. Secondly, the numerical evidence suggests that for fixed $s$ and varying $\frac{L}{\ell}$, the influence of the magnitude of $\frac{L}{\ell}$ on $\alpha^\star_\gamma$ is negligible. It is straightforward to show this analytically, see the simplified versions of (\ref{Da0_Lin}) and (\ref{Da1_Lin}).

\begin{figure}[!ht]
\begin{center}
\includegraphics[width=1.0\textwidth]{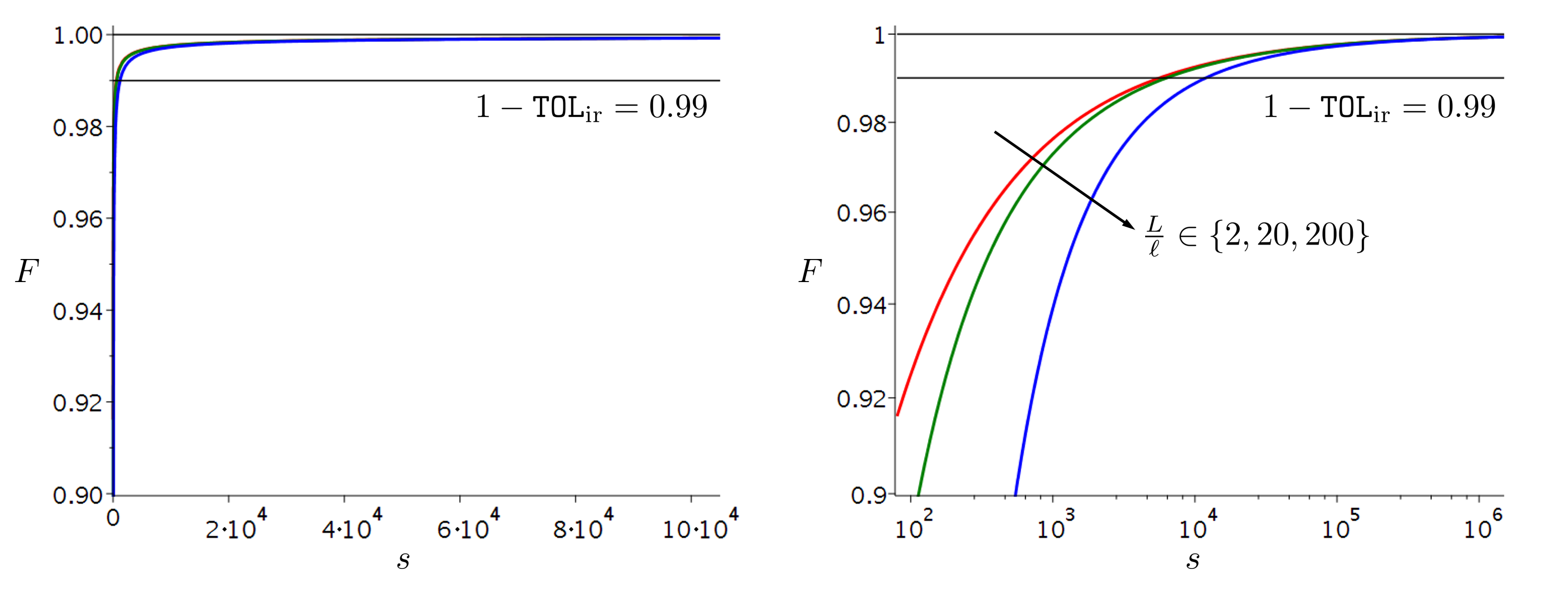}
\end{center}
\caption{Function $F$ in (\ref{Fs_Lin}) plotted on the interval $s\in[10^2,10^6]$ for the three different ratios $\frac{L}{\ell}\in\{2,20,200\}$ in the normal (on the left) and logarithmic (on the right) scales; the horizontal line $1-\mathtt{TOL}_\mathrm{ir}$, with $\mathtt{TOL}_\mathrm{ir}:=0.01$, represents the irreversibility threshold for obtaining the optimal penalty parameter $s_\mathrm{opt}$.}
\label{Fig5}
\end{figure}

We are now in the position to compute the quantity $G_c^{-1}\widetilde{E}_S(\alpha^\star_\gamma)$. Substituting (\ref{aStarGamm_Lin}) in $\widetilde{E}_S$ in (\ref{Es1d_gam}), we arrive at 
\begin{align*}
G_c^{-1}\widetilde{E}_S(\alpha^\star_\gamma)
=&1-\frac{3}{16s}\left( \frac{L}{\ell}-2 \right)
-\frac{3}{4s}\exp(2\sqrt{s}) \frac{ \exp\left( 2\sqrt{s}\left( \frac{L}{\ell}-2 \right) \right)-1 }{ \exp\left( 2\sqrt{s}\frac{L}{\ell} \right)-1 } \\
-&\frac{3}{4s^\frac{3}{2}} \left( s-\frac{1}{8} \right) \frac{ \exp\left( 2\sqrt{s}\frac{L}{\ell} \right)+1 }{ \exp\left( 2\sqrt{s}\frac{L}{\ell} \right)-1 }
-\frac{3}{32s^\frac{3}{2}} \frac{ \exp\left( 2\sqrt{s}\left( \frac{L}{\ell}-2 \right) \right)+\exp(4\sqrt{s}) }{ \exp\left( 2\sqrt{s}\frac{L}{\ell} \right)-1 }.
\end{align*}
The obtained expression is a function of both $s$ and $\frac{L}{\ell}$. However, assuming that $s\gg1$ and $\frac{L}{\ell}\geq2$, the above can be simplified significantly. As a result, we obtain function $F$ defined as follows
\begin{equation}
G_c^{-1}\widetilde{E}_S(\alpha^\star_\gamma)\approx F(s):=
\left\{
\begin{tabular}{ll}
$\displaystyle 1-\frac{3}{4s^\frac{1}{2}}$,  &  if $\displaystyle \frac{L}{\ell}=2$, \\ [0.2cm]
$\displaystyle 1-\frac{3}{16s}\left( \frac{L}{\ell}-2 \right)
-\frac{3}{4s^\frac{3}{2}} \left( s-\frac{1}{8} \right)$,  &  if $\displaystyle \frac{L}{\ell}>2$.
\end{tabular}
\right.
\label{Fs_Lin}
\end{equation}
The plots of $F$ in (\ref{Fs_Lin}) as a function of $s$ in different scales are depicted in Figure \ref{Fig5}. 

\begin{figure}[!ht]
\begin{center}
\includegraphics[width=1.0\textwidth]{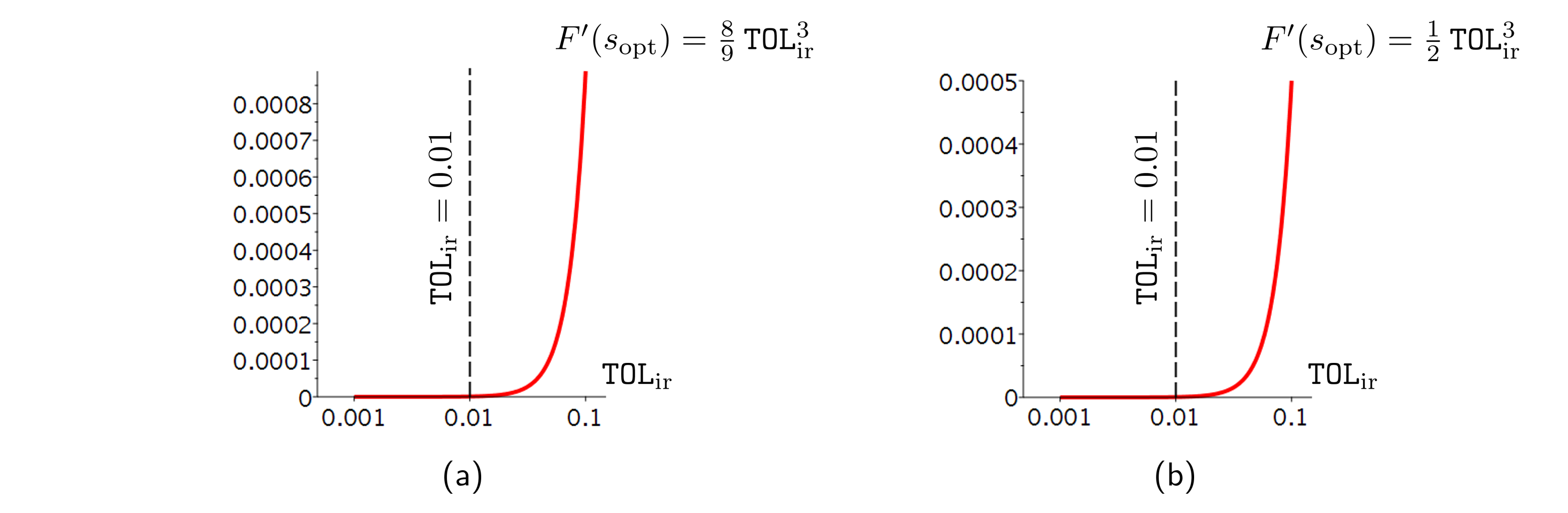}
\end{center}
\caption{The plot of $F^\prime(s_\mathrm{opt})$ as a function of $\mathtt{TOL}_\mathrm{ir}\in[0.001,0.1]$ for: (a) the linear model, equation (\ref{dFsOpt_Lin}), and (b) the quadratic model, equation (\ref{dFsOpt_Sq}); in either case $\mathtt{TOL}_\mathrm{ir}:=0.01$ can be treated as a practical irreversibility threshold.}
\label{Fig6}
\end{figure}

The following observations can be made:
\begin{itemize}
	\item[1.] The increase of $s$ yields a monotonic increase of $F$. In particular, when $s\rightarrow +\infty$ (and, hence, $\gamma\rightarrow +\infty$), $F$ asymptotically approaches $1$. This mimics the $\Gamma$-convergence result for the corresponding minimizer $\alpha^\star_\gamma$, which is given in terms of the penalty parameter $s$.
	\item[2.] We call the intersection point between $y=F(s)$ and $y=1-\mathtt{TOL}_\mathrm{ir}$, where $0<\mathtt{TOL}_\mathrm{ir}\ll1$, the optimal value of $s$ and denote it as $s_\mathrm{opt}$. Herein, we also call $\mathtt{TOL}_\mathrm{ir}$ the (user-prescribed) irreversibility tolerance threshold. That is, technically, $s_\mathrm{opt}$ is the solution of the equation $F(s)=1-\mathtt{TOL}_\mathrm{ir}$ and conceptually, it represents the penalty parameter sufficient to provide for the solution $\alpha^\star_\gamma$ the recovery of the desired reference $\Gamma$-convergence value $1$ with the error of $\mathtt{TOL}_\mathrm{ir}\cdot 100\%$.
	\item[3.] In the case $\frac{L}{\ell}>2$, function $F$ depends on the ratio $\frac{L}{\ell}$, that is, the value $s_\mathrm{opt}$ we are interested in is size dependent. The plots depicted in Figure \ref{Fig5}, however, suggest that when a practical range of $\frac{L}{\ell}$ is considered, $s_\mathrm{opt}$ may be viewed as a $\frac{L}{\ell}$-independent result. Indeed, for the irreversibility tolerance $\mathtt{TOL}_\mathrm{ir}:=0.01$, the values $s_\mathrm{opt}$ in the two cases of $\frac{L}{\ell}=2$ and e.g.\ $\frac{L}{\ell}=200$ are quite close, as can be observed in the left (normal scale) plot of $F$ in Figure \ref{Fig5}. The right (log scale) plot in Figure \ref{Fig5} shows that they are of the same order of magnitude, the values being $s_\mathrm{opt}=5625$ and $s_\mathrm{opt}\approx11890$, respectively.
	\item[4.] If we agree to take $F$ in (\ref{Fs_Lin}) with $\frac{L}{\ell}=2$ as the reference case for determining $s_\mathrm{opt}$ (i.e.\ the dependence of $s_\mathrm{opt}$ on $\frac{L}{\ell}>2$ is neglected), we can argue that $\mathtt{TOL}_\mathrm{ir}:=0.01$ may be viewed as a practical reference irreversibility threshold. Indeed, already the left plot of Figure \ref{Fig5} heuristically suggests that smaller values of $\mathtt{TOL}_\mathrm{ir}$ do not improve the $\Gamma$-convergence error significantly. Rigorously, this can be shown as follows. The solution of the equation $F(s)=1-\mathtt{TOL}_\mathrm{ir}$ yields
\begin{equation}
s_\mathrm{opt}:=\frac{9}{16\,\mathtt{TOL}_\mathrm{ir}^2},
\label{sOpt_Lin}
\end{equation}
where $0<\mathtt{TOL}_\mathrm{ir}\ll1$. Then, inserting (\ref{sOpt_Lin}) into $F^\prime(s)=\frac{3}{8}s^{-\frac{3}{2}}$, one gets
\begin{equation}
F^\prime(s_\mathrm{opt}):=\frac{8}{9}\mathtt{TOL}_\mathrm{ir}^3,
\label{dFsOpt_Lin}
\end{equation}
whose plot is depicted in Figure \ref{Fig6}(a). It can be seen that in the range of $\mathtt{TOL}_\mathrm{ir}\in[0.001,0.1]$, neither $F^\prime$, nor its slope change for all $\mathtt{TOL}_\mathrm{ir}$ less than the proposed reference value $0.01$. In section \ref{SecNum}, the numerical evidence of this assumption will also be illustrated.
\end{itemize} 

\subsubsection{Quadratic model, solution to (\ref{bvp_s})}
\label{SqMod_gam}
The exact solution to problem (\ref{bvp_s}) can be constructed using the Newton-Raphson procedure in only one iteration, see appendix \ref{Ap1_SqMod}, thus reading 
\begin{equation}
\alpha^\star_\gamma(x)=\alpha^\star(x)+\Delta\alpha_0(x),
\label{aStarGamm_Sq}
\end{equation}
where $\Delta\alpha_0$ is explicitly given in appendix \ref{Ap1_SqMod} by equation (\ref{Da0_Sq}).

\begin{figure}[!ht]
\begin{center}
\includegraphics[width=1.0\textwidth]{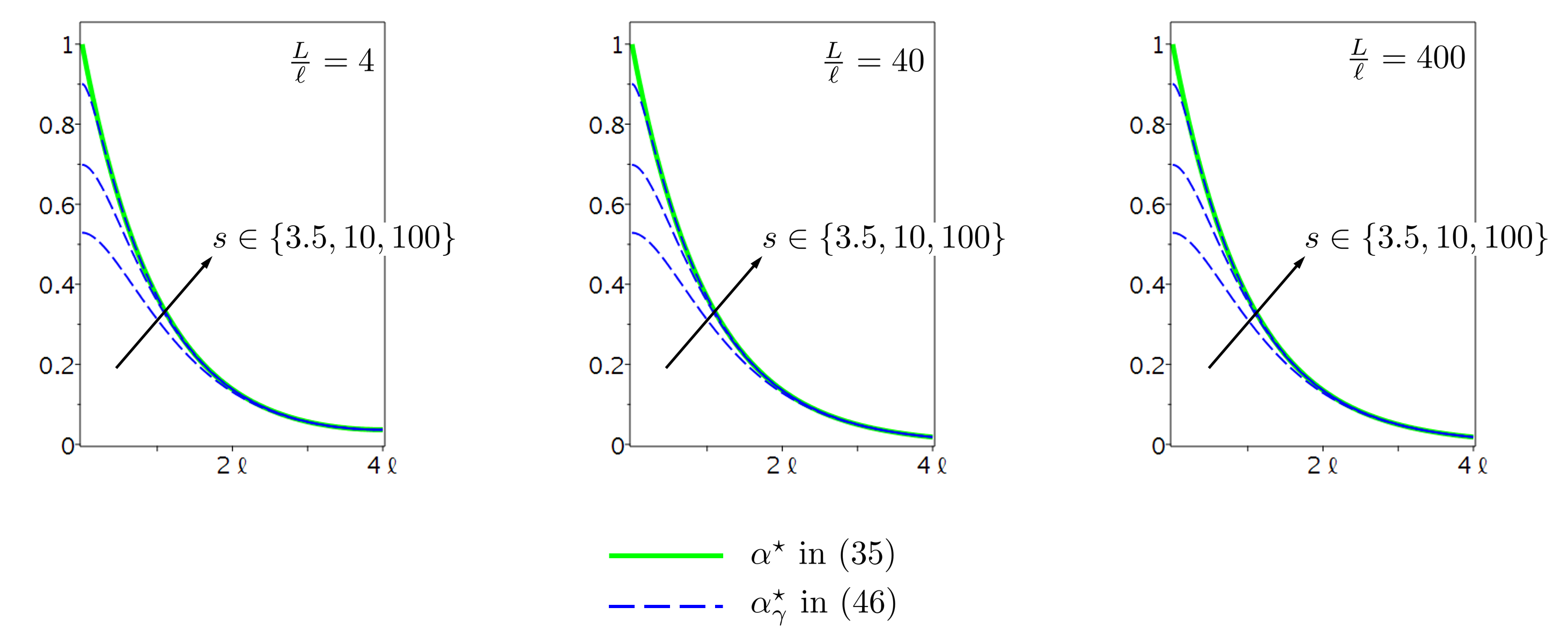}
\end{center}
\caption{The functions $\alpha^\star$ in (\ref{aStar1d_Sq}) and $\alpha^\star_\gamma$ in (\ref{aStarGamm_Sq}) plotted on the interval $x\in[0,4\ell]$ for $\frac{L}{\ell}\in\{4,40,400\}$; for fixed $\frac{L}{\ell}$, $\alpha^\star_\gamma$ is depicted for different values of the penalty parameter $s>0$.}
\label{Fig7}
\end{figure}

In Figure \ref{Fig7}, we compare the constructed profile $\alpha^\star_\gamma$ in (\ref{aStarGamm_Sq}) and the reference $\alpha^\star$ in (\ref{aStar1d_Sq}) for three different magnitudes of $\frac{L}{\ell}>2$. Also, for fixed $\frac{L}{\ell}$, $\alpha^\star_\gamma$ is depicted depending on the penalty parameter $s>0$. For the considered cases $\frac{L}{\ell}\in\{4,40,400\}$, both $\alpha^\star$ and $\alpha^\star_\gamma$ are defined, respectively, on the intervals $x\in[0,4\ell]$, $[0,40\ell]$ and $[0,400\ell]$, yet for comparison purposes all plots are restricted to the interval $[0,4\ell]$. Similar to the linear model case, we conclude here that for a fixed $\frac{L}{\ell}$, the desired point-wise convergence of $\alpha^\star_\gamma$ to $\alpha^\star$ when $s$ increases is demonstrated. It is also numerically evident that for fixed $s$ and varying $\frac{L}{\ell}$, the dependence of $\alpha^\star_\gamma$ on $\frac{L}{\ell}$ is negligible. This can be also confirmed analytically, see the simplified representation for (\ref{Da0_Sq}).

\begin{figure}[!ht]
\begin{center}
\includegraphics[width=1.0\textwidth]{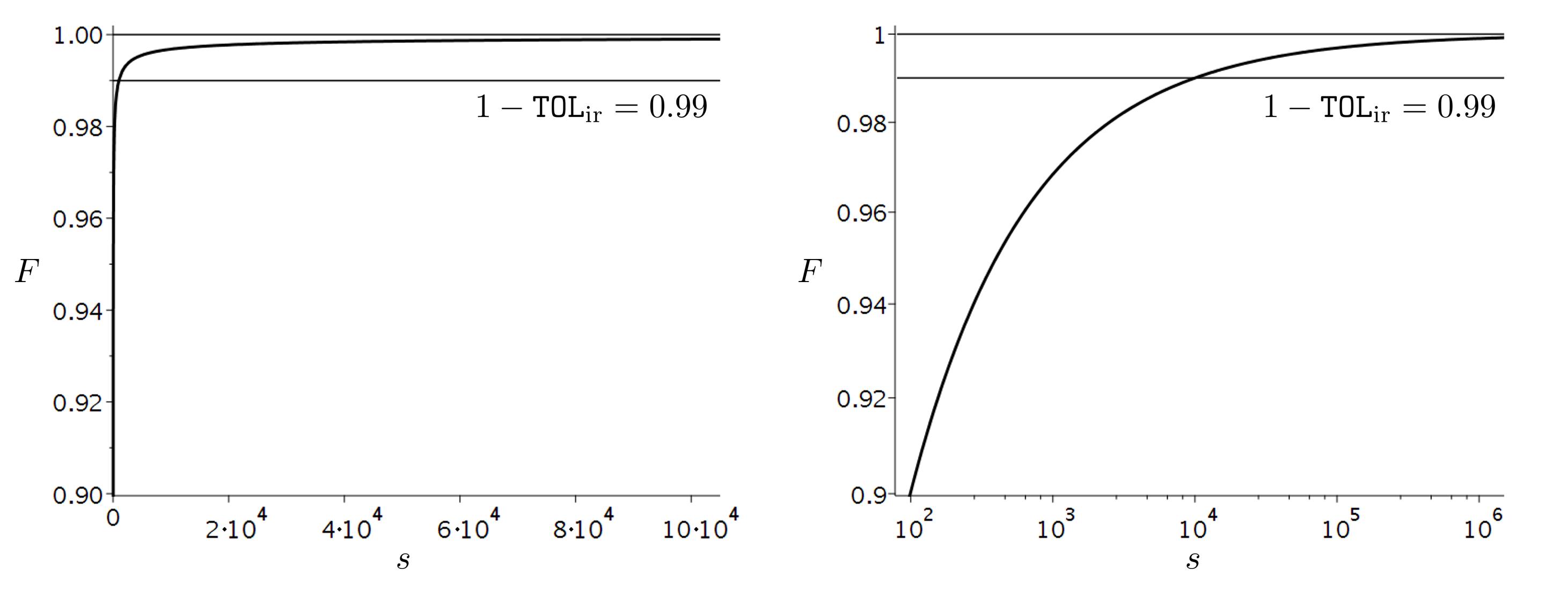}
\end{center}
\caption{Function $F$ in (\ref{Fs_Sq}) plotted on the interval $s\in[10^2,10^6]$ in the normal (on the left) and logarithmic (on the right) scales; the horizontal line $1-\mathtt{TOL}_\mathrm{ir}$, with $\mathtt{TOL}_\mathrm{ir}:=0.01$, represents the irreversibility threshold for obtaining the optimal penalty parameter $s_\mathrm{opt}$.}
\label{Fig8}
\end{figure}

Substituting (\ref{aStarGamm_Sq}) in $\widetilde{E}_S$ in (\ref{Es1d_gam}), we arrive at 
\begin{equation*}
G_c^{-1}\widetilde{E}_S(\alpha^\star_\gamma)
=\frac{ \exp\left( 2\frac{L}{\ell} \right)-1 }{ \exp\left( 2\frac{L}{\ell} \right)+1 }
-\frac{1}{(s+1)^\frac{1}{2}} 
 \left( \frac{ \exp\left( 2\frac{L}{\ell} \right)-1 }{ \exp\left( 2\frac{L}{\ell} \right)+1 } \right)^2
 \frac{ \exp\left( 2\sqrt{s+1}\frac{L}{\ell} \right)+1 }{ \exp\left( 2\sqrt{s+1}\frac{L}{\ell} \right)-1 }.
\end{equation*}
The obtained expression is a function of both $s$ and $\frac{L}{\ell}$. Assuming again that $s\gg1$ and $\frac{L}{\ell}\geq2$, the above can be simplified to finally obtain the function $F$,
\begin{equation}
G_c^{-1}\widetilde{E}_S(\alpha^\star_\gamma)\approx F(s):=
1-\frac{1}{(s+1)^\frac{1}{2}}. 
\label{Fs_Sq}
\end{equation}
The plots of $F$ in (\ref{Fs_Sq}) as a function of $s$ in the different scales are depicted in Figure \ref{Fig8}.

All observations in section \ref{LinMod_gam} for $F$ in the linear case are valid for $F$ in (\ref{Fs_Sq}), except that the present $F$ is independent on the ratio $\frac{L}{\ell}\geq2$. The solution of the equation $F(s)=1-\mathtt{TOL}_\mathrm{ir}$ reads
\begin{equation}
s_\mathrm{opt}:=\frac{1}{\mathtt{TOL}_\mathrm{ir}^2}-1,
\label{sOpt_Sq}
\end{equation}
such that 
\begin{equation}
F^\prime(s_\mathrm{opt}):=\frac{1}{2}\mathtt{TOL}_\mathrm{ir}^3,
\label{dFsOpt_Sq}
\end{equation}
see Figure \ref{Fig6}(b). Similarly to the linear case -- see equation (\ref{dFsOpt_Lin}) and the corresponding plot in Figure \ref{Fig6}(a) -- $\mathtt{TOL}_\mathrm{ir}:=0.01$ may be taken as a practical reference irreversibility threshold for the quadratic model.

\subsubsection{Final remarks}
Once the problem-independent $s_\mathrm{opt}$ in (\ref{sOpt_Lin}) and (\ref{sOpt_Sq}) is given, the relation (\ref{s_variable}) is used for computing the actual penalty parameter $\gamma$ for the penalized formulation in (\ref{argmin1}). It is defined as
\begin{equation}
\gamma_\mathrm{opt}:=\frac{2}{c_{\sf w}}\frac{G_c}{\ell}s_\mathrm{opt}
=\left\{
\begin{tabular}{ll}
$\displaystyle \frac{G_c}{\ell}\frac{27}{64\,\mathtt{TOL}_\mathrm{ir}^2}$, & when ${\sf w}(\alpha)=\alpha$,   \\  [0.3cm]
$\displaystyle \frac{G_c}{\ell}\left( \frac{1}{\mathtt{TOL}_\mathrm{ir}^2}-1 \right)$,  & when ${\sf w}(\alpha)=\alpha^2$,
\end{tabular}
\right.
\label{gamOpt}
\end{equation}
with $0<\mathtt{TOL}_\mathrm{ir}\ll1$ as user-prescribed irreversibility tolerance threshold. 

Recall that one of the assumptions behind the derivation idea for $s_\mathrm{opt}$ is that $\alpha_{n-1}$ already represents the fully developed crack. Therefore, the obtained $s_\mathrm{opt}$ may be insufficient to enforce irreversibility of $\alpha$ when $\alpha_{n-1}<1$. It is also worth noticing that the obtained result can be extrapolated to problems in higher dimensions. In section \ref{SecNum}, we will bring the numerical evidence of this.

\subsection{$\rho$-penalization: optimal profile for the linear model}
\label{Penal2}
In this section, we construct the minimizer $\alpha^\star_\rho$ of the one-dimensional counterpart of $\widehat{E}_S$ in (\ref{Es_rho}) reading
\begin{equation}
\widehat{E}_S(\alpha):=
\frac{3}{4}G_c \int_0^L
\left(\frac{\alpha}{\ell}+\ell \left( \alpha^\prime \right)^2\right)\mathrm{d}x
+\rho\int_0^L\langle \alpha \rangle_{-}^2 \; \mathrm{d}x,
\label{Es1d_rho}
\end{equation}
(note again that owing to symmetry, we consider the half-length representation). We then establish the related $\Gamma$-convergence property and derive the lower bound for the penalty constant $\rho\gg1$. 

The weak formulation reads: find $\alpha\in\{H^1(0,L): \alpha(0)=1\}$ such that
\begin{equation}
\widetilde{E}_S^\prime(\alpha,\beta):=\frac{3}{4}G_c
\int_0^L \left(
\frac{1}{\ell}\beta+2\ell\alpha^\prime \beta^\prime
\right)\mathrm{d}x
+2\rho\int_0^L\langle \alpha \rangle_{-}\beta\;\mathrm{d}x
\overset{!}{=}0, \quad \forall \beta\in H^1(0,L): \beta(0)=0.
\label{weak_rho}
\end{equation}
The resulting boundary value problem for $\alpha^\star_\rho$ is as follows:
\begin{equation}
\left\{
\begin{tabular}{rl}
$\displaystyle -\ell^2\alpha^{\prime\prime}+\frac{1}{2}
+\frac{4}{3}\frac{\rho\ell}{G_c} \langle \alpha \rangle_{-}=0$ & in $(0,L)$,   \\ 
$\alpha=1$ & at $x=0$, \\ 
$\alpha^\prime=0$ & at $x=L$.
\end{tabular}
\right.
\label{bvp_rho}
\end{equation}
As before, we introduce the new (dimensionless) penalty parameter:
\begin{equation}
\boxed
{
\frac{4}{3}\frac{\rho\ell}{G_c}=:r,
}
\label{r_variable}
\end{equation}
such that $r>0$, and rewrite (\ref{bvp_rho}) compactly as follows: 
\begin{equation}
\left\{
\begin{tabular}{rl}
$\displaystyle -\ell^2\alpha^{\prime\prime}+\frac{1}{2}
+r \langle \alpha \rangle_{-}=0$ & in $(0,L)$,   \\ 
$\alpha=1$ & at $x=0$, \\ 
$\alpha^\prime=0$ & at $x=L$.
\end{tabular}
\right.
\label{bvp_r}
\end{equation}

\begin{figure}[!ht]
\begin{center}
\includegraphics[width=1.0\textwidth]{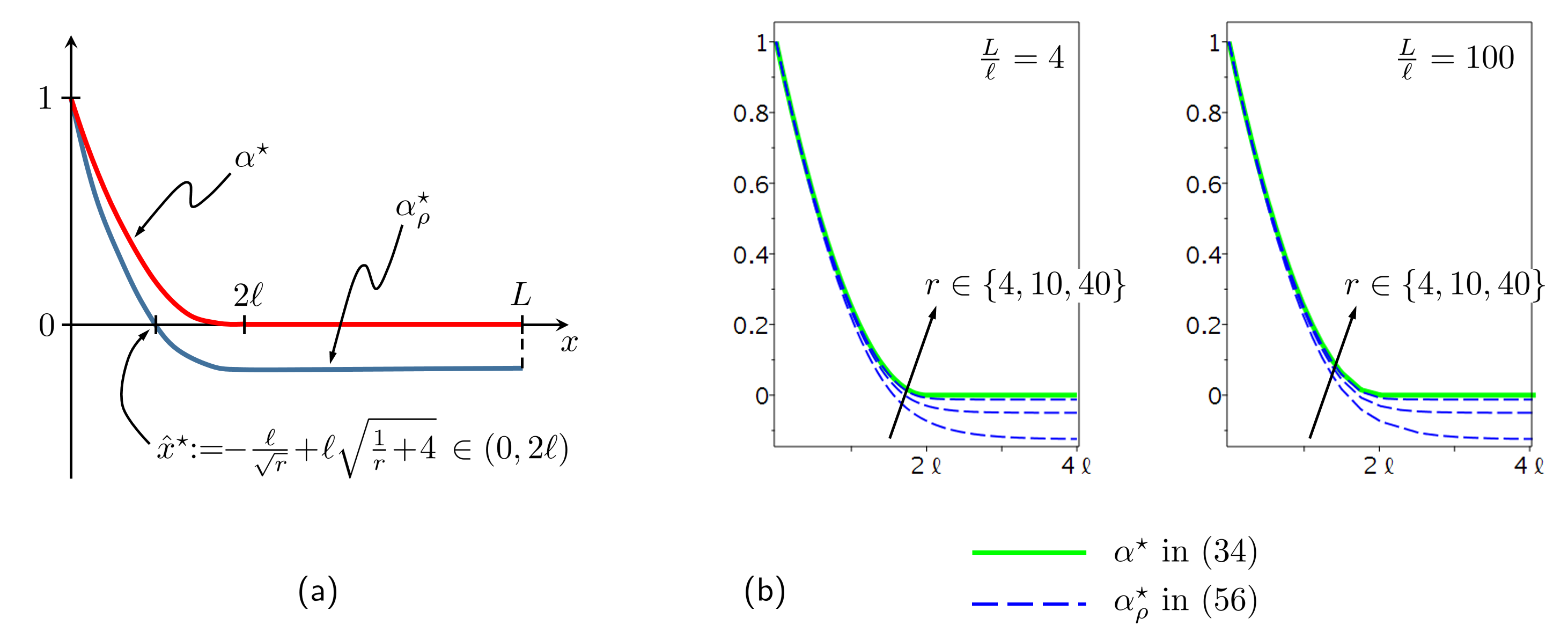}
\end{center}
\caption{(a) Sketch of the reference profile $\alpha^\star$ in (\ref{aStar1d_Lin}), $\alpha^\star_\rho$ in (\ref{aStarRho}) and the corresponding limiting point $\hat{x}^\star$: when $\hat{x}^\star$ converges to $2\ell$, the point-wise convergence of $\alpha^\star_\rho$ to $\alpha^\star$ is expected; (b) The numerical evidence of the convergence behavior of $\alpha^\star_\rho$ to $\alpha^\star$ in terms of the penalty parameter $r$.}
\label{Fig9}
\end{figure}

The (exact) solution of the semi-linear problem in (\ref{bvp_r}) is constructed in Appendix \ref{Ap2}. It reads
\begin{equation}
\alpha^\star_\rho(x):=\left\{
\begin{tabular}{ll}
$\displaystyle \frac{1}{4\ell^2}x^2+R(\hat{x}^\star)x+1$, & in $(0,\hat{x}^\star)$,   \\ [0.2cm]
$\displaystyle T(\hat{x}^\star)\exp\left(-\sqrt{r}\frac{x}{\ell}\right) -\frac{1}{2r}$, & in $(\hat{x}^\star,L)$,
\end{tabular}
\right.
\label{aStarRho}
\end{equation}
where
\begin{equation*}
\hat{x}^\star:=-\frac{\ell}{\sqrt{r}}+\ell\sqrt{\frac{1}{r}+4}
\end{equation*}
and
\begin{equation*}
R(\hat{x}^\star):=-\frac{1}{2\ell}\sqrt{\frac{1}{r}+4},
\end{equation*}
\begin{equation*}
T(\hat{x}^\star):=\frac{1}{2r}\exp\left( -1+\sqrt{1+4r} \right).
\end{equation*}
Figure \ref{Fig9}(a) sketches the link between the reference optimal profile $\alpha^\star$ in (\ref{aStar1d_Lin}) and the constructed profile $\alpha^\star_\rho$ in (\ref{aStarRho}). More precisely, it highlights the anticipation that the convergence of $\hat{x}^\star$ to $2\ell$ when $r\rightarrow+\infty$ must result in the point-wise convergence of $\alpha^\star_\rho$ to the corresponding $\alpha^\star$. Figure \ref{Fig9}(b), where for fixed ratios $\frac{L}{\ell}\in\{4,100\}$, $\alpha^\star_\rho$ is plotted for different values of parameter $r$, presents the related numerical evidence.  

\begin{figure}[!ht]
\begin{center}
\includegraphics[width=1.0\textwidth]{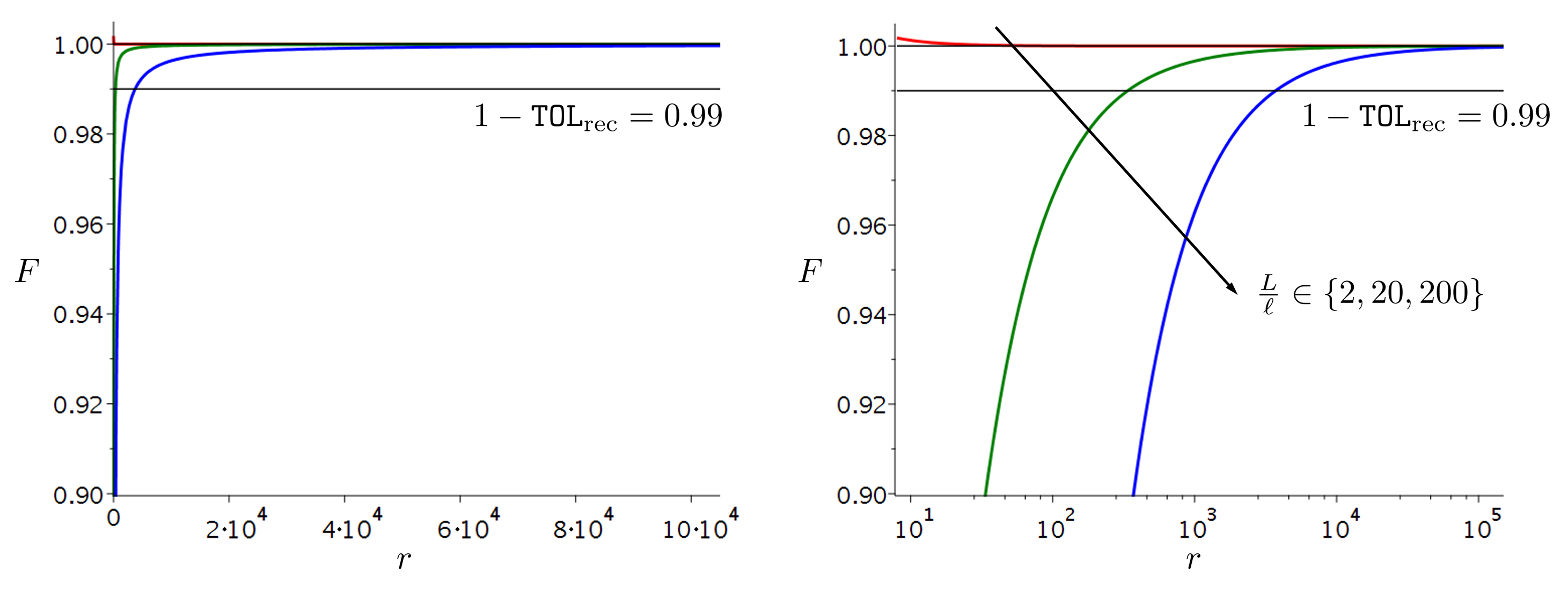}
\end{center}
\caption{Function $F$ in (\ref{Fr}) plotted on the interval $r\in[10^1,10^5]$ for $\frac{L}{\ell}\in\{2,20,200\}$ in the normal (on the left) and logarithmic (on the right) scales; the horizontal line $1-\mathtt{TOL}_\mathrm{rec}$, with $\mathtt{TOL}_\mathrm{rec}:=0.01$, represents the recovery threshold for obtaining the optimal penalty parameter $r_\mathrm{opt}$.}
\label{Fig10}
\end{figure}

The final step of the presented analysis is the derivation of the optimal lower bound for $r$ (and, hence, for $\rho$). As before, we insert (\ref{aStarRho}) into $\widehat{E}_S(\alpha)$ in (\ref{Es1d_rho}) and arrive at
\begin{align*}
G_c^{-1}\widehat{E}_S(\alpha^\star_\rho)=&
\left(1+\frac{1}{4r}\right)^\frac{3}{2}-\frac{3}{16r}\left(\frac{L}{\ell}\right)
+\frac{1}{16r^\frac{3}{2}}-\frac{3}{16r^\frac{3}{2}}
\exp\left( 2\sqrt{1+4r}-2\sqrt{r}\frac{L}{\ell}-2 \right) \nonumber \\
\approx& \left(1+\frac{1}{4r}\right)^\frac{3}{2}-\frac{3}{16r}\left(\frac{L}{\ell}\right)
+\frac{1}{16r^\frac{3}{2}}.
\end{align*}
The above simplification holds true, assuming $r$ is large and $\frac{L}{\ell}\geq2$. Further simplification applies using the Taylor expansion of the term $\left(1+\frac{1}{4r}\right)^\frac{3}{2}=1+\frac{3}{8r}+O(\frac{1}{r^2})$, such that we eventually have
\begin{equation}
G_c^{-1}\widehat{E}_S(\alpha^\star_\rho)
\approx F(r):=1-\frac{3}{16r}\left(\frac{L}{\ell}-2\right)+\frac{1}{16r^\frac{3}{2}}.
\label{Fr}
\end{equation}
The plots of $F$ in (\ref{Fr}) as a function of $r$ in the different scales are depicted in Figure \ref{Fig10}. 

\begin{figure}[!ht]
\begin{center}
\includegraphics[width=1.0\textwidth]{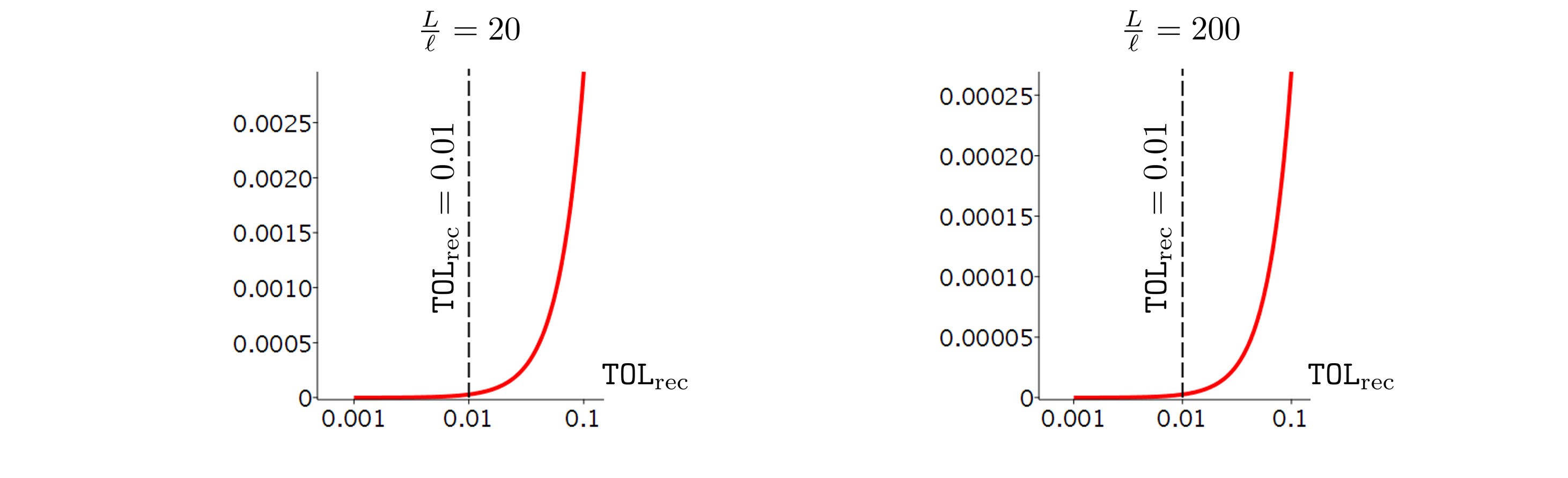}
\end{center}
\caption{The plot of $F^\prime_a(r_\mathrm{opt})$ as a function of $\mathtt{TOL}_\mathrm{rec}\in[0.001,0.1]$ for $\frac{L}{\ell}\in\{20,200\}$; this is to argue that in either case $\mathtt{TOL}_\mathrm{rec}\leq0.01$ can be treated as the reference recovery threshold.}
\label{Fig11}
\end{figure}

The observations that lead to the calculation of the optimal $r$ (and, hence, $\rho$) are as follows:
\begin{itemize}
	\item[1.] Function $F$ asymptotically approaches $1$ when $r\rightarrow+\infty$, what mimics the desired $\Gamma$-convergence property of the constructed minimizer $\alpha^\star_\rho$. In particular, $1$ is the lower limit for $F$ when $\frac{L}{\ell}=2$ and the upper limit for all $\frac{L}{\ell}>2$. This can straightforwardly be seen from (\ref{Fr}). 
	\item[2.] For $\frac{L}{\ell}>2$, the optimal value of $r$, denoted as $r_\mathrm{opt}$, is a solution of the equation $F(r)=1-\mathtt{TOL}_\mathrm{rec}$, where $0<\mathtt{TOL}_\mathrm{rec}\ll1$ stands for the user-prescribed recovery tolerance threshold. In this context, $r_\mathrm{opt}$ represents the penalty parameter sufficient to obtain for the solution $\alpha^\star_\rho$ the recovery of the desired reference $\Gamma$-convergence value $1$ with the error of $\mathtt{TOL}_\mathrm{rec}\cdot 100\%$.
	\item[3.] Since $F$ depends on the ratio $\frac{L}{\ell}$, $r_\mathrm{opt}$ depends on the structural dimension $L$: see the right plot in Figure \ref{Fig10}, where a large difference in the order of magnitude for $r_\mathrm{opt}$ is clearly visible in the two cases of $\frac{L}{\ell}=20$ and $200$. That is, in contrast to the $\gamma$-penalization case, the dependence of the results on $\frac{L}{\ell}$ cannot be neglected.  
	\item[4.] The equation $F(r)=1-\mathtt{TOL}_\mathrm{rec}$ for obtaining $r_\mathrm{opt}$ is a cubic polynomial equation in $y:=r^\frac{1}{2}$. Its analytical solution is available but very cumbersome. Without loss of generality, we consider the approximation $F_a$ of $F$ in (\ref{Fr}) obtained by neglecting the higher order term $\frac{1}{16r^\frac{3}{2}}$ (this is justified by the assumption that $r$ is large). The equation $F_a(r)=1-\mathtt{TOL}_\mathrm{rec}$ is linear in $y$, and the corresponding solution reads
\begin{equation}
r_\mathrm{opt}:=\frac{3}{16}\frac{\frac{L}{\ell}-2}{\mathtt{TOL}_\mathrm{rec}}.
\label{rOpt_Lin}
\end{equation}
It is straightforward to check that (\ref{rOpt_Lin}) is a sufficiently good approximation of the solution to the original cubic equation for large $r$.  
 	\item[5.] As a final step, let us argue that any $\mathtt{TOL}_\mathrm{rec}\leq0.01$ may be taken as a practical recovery threshold. For this, we substitute (\ref{rOpt_Lin}) in $F^\prime_a(r)=\frac{3}{16r^2}\left( \frac{L}{\ell}-2\right)$, thus yielding
\begin{equation}
F^\prime_a(r_\mathrm{opt}):=\frac{16}{3\left( \frac{L}{\ell}-2 \right)}
\,\mathtt{TOL}_\mathrm{rec}^2,
\label{dFrOpt}
\end{equation}
and plot $F^\prime_a$ in the range of $\mathtt{TOL}_\mathrm{rec}\in[0.001,0.1]$ for two different values of $\frac{L}{\ell}$, see Figure \ref{Fig11}. It can be seen that for $\frac{L}{\ell}\in\{20,200\}$, both $F^\prime_a$ and its slope have a minor change for all $\mathtt{TOL}_\mathrm{rec}$ less than $0.01$. The numerical evidence will be illustrated for the corresponding examples in section \ref{SecNum}.
\end{itemize}

With $r_\mathrm{opt}$ given by (\ref{rOpt_Lin}), the relation (\ref{r_variable}) is used for computing the actual penalty parameter $\rho$ for the corresponding penalized recovery formulation:
\begin{equation}
\rho_\mathrm{opt}:=\frac{3}{4}\frac{G_c}{\ell}r_\mathrm{opt}
=\frac{G_c}{\ell} \frac{9\left( \frac{L}{\ell}-2 \right)}{64\,\mathtt{TOL}_\mathrm{rec}}.
\label{rhoOpt}
\end{equation}
As mentioned for the $\gamma$-penalization result, we extrapolate the obtained results to higher dimensions. In this case, the domain size parameter $L$ is replaced by either $\mathrm{diam}(\Omega)$ or $L_{\max}$.

\section{Numerical studies}
\label{SecNum}
The obtained results for $\gamma_\mathrm{opt}$ and $\rho_\mathrm{opt}$ are extrapolated to problems in two dimensions. To verify these findings the following numerical experiments are considered: the so-called single edge notched (SEN) specimen under {\em shear}, and the Sneddon-Lowengrub problem. 

The first one is a crack propagation problem under external displacement-controlled loading, where the pre-existing crack is modeled discretely and the propagating crack is represented by the phase-field evolution. The problem setup is simple, but the failure pattern is not symmetric, being the result of a non-trivial combination of local tension-compression within the specimen during shear. The loading-unloading and pure loading regimes are simulated to fully test the $\gamma$-penalized irreversibility.

The Sneddon-Lowengrub problem in its original discrete setting is a crack {\em opening} problem, that is, a problem where the internal pressure applied to the pre-existing crack faces only opens the crack but does not result in its propagation. The problem is of interest and importance by the following aspects. First, since the phase-field treatment implies modeling of the pre-existing crack by means of the phase-field, we are able to test our $\rho$-penalized recovery procedure. Secondly, the original discrete problem has an exact analytical solution, also for the so-called crack opening displacement (COD). The phase-field representation of the COD involves both the displacement field and the crack phase-field. Therefore, the comparison of the results enables to explicitly assess the accuracy of the phase-field formulation.

Apart from testing the aforementioned $\gamma$-penalized irrreversibility and $\rho$-penalized recovery, for both problems the qualitative and quantitative impact of the ingredients of $\mathcal{E}$ in (\ref{RegVAF1}) such as the tension-compression split and the local part of the fracture energy density will be evaluated. Finally, for the SEN shear test we will compare our penalized irreversibility results, that is, the solution of (\ref{Weak1}) with the solution of (\ref{Weak4}) where the history field of Miehe and co-workers is used. 

We employ the numerical package FreeFem++ \cite{FreeFem}. Both the displacement field $\bm u$ and phase-field $\alpha$ are approximated using $P_1$-triangles. The error tolerance for the staggered and the (nested) Newton-Raphson iterative solution processes are prescribed as $\mathtt{TOL}_\mathrm{Stag}:=10^{-4}$ and $\mathtt{TOL}_\mathrm{NR}:=10^{-6}$, respectively, see section \ref{StagSol} for details.

\subsection{SEN specimen under shear}
\begin{figure}[!ht]
\begin{center}
\includegraphics[width=1.0\textwidth]{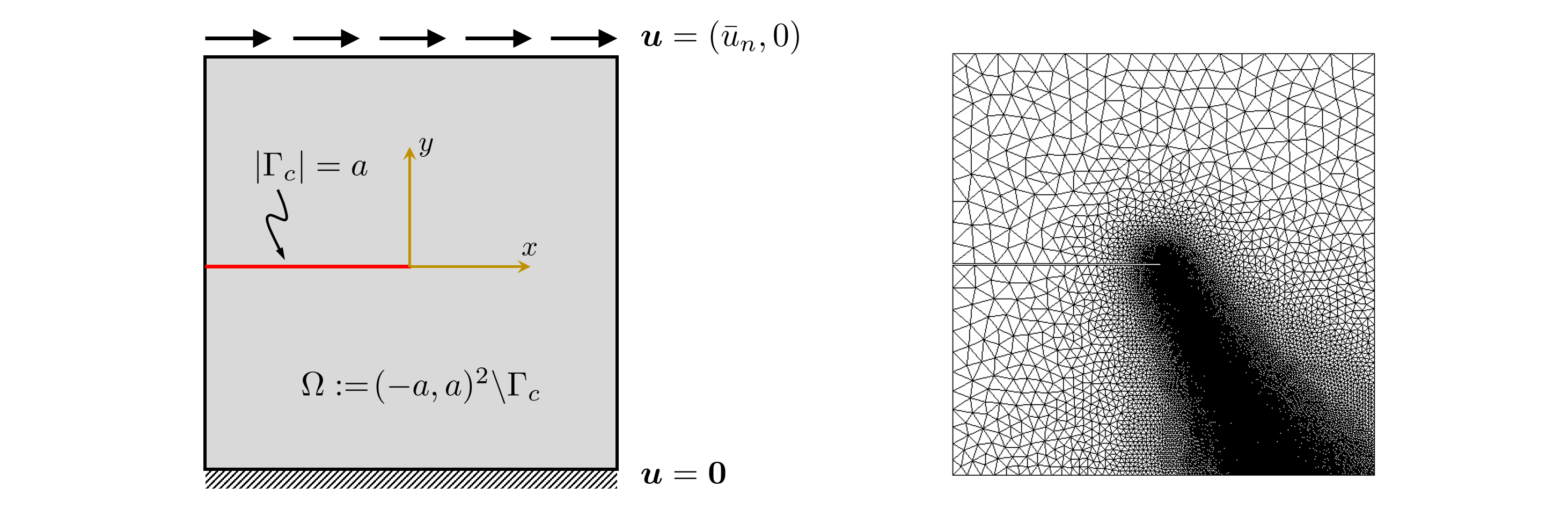}
\end{center}
\caption{Geometry, loading setup and the (pre-adapted) finite element mesh for the single edge notched (SEN) specimen subject to shear.}
\label{Fig12}
\end{figure}

The problem setup is depicted in the left plot of Figure \ref{Fig12}. This benchmark example was originally considered in \cite{Bourdin2000}, and later adopted with minor modifications in many related papers, see e.g.\ \cite{Miehe2010b,Borden2014,Ambati2015review,Gerasimov2016}.  

The specimen geometry, material and loading data are taken as in \cite{Miehe2010b} and are as follows: $a=1$ mm, $E=210$ GPa, $\nu=0.3$, $G_c=2700$ N/m. We set the length scale parameter as $\ell=\frac{2a}{200}=0.01$ mm. Plane strain is assumed. A horizontal displacement loading $\bar{u}_n$ (which monotonically increases in a pure loading regime, and increases-decreases in a loading-unloading regime) is prescribed on the upper specimen edge, with the vertical component being restrained. The bottom edge is completely fixed. As the crack propagation pattern is known from earlier phase-field studies, our finite element mesh is pre-adapted in the region where the phase-field evolution is expected, see the right plot in Figure \ref{Fig12}. The mesh contains 33193 elements and 16761 nodes. The minimal mesh size in the localization zone of $\alpha$ is $h\approx 0.0025$ mm, that is, $h\ll\ell$ is fulfilled.

\subsubsection{Parametric studies for $\gamma$}
As argued in section \ref{Penal1}, we take $\mathtt{TOL}_\mathrm{ir}:=0.01$. Then, recalling (\ref{gamOpt}), the penalty parameter $\gamma$ which we term optimal is obtained:
\begin{equation*}
\gamma_\mathrm{opt}\approx
\left\{
\begin{tabular}{ll}
$1.14\cdot 10^{12}$, & when ${\sf w}(\alpha)=\alpha$,   \\
$2.7\cdot 10^{12}$,  & when ${\sf w}(\alpha)=\alpha^2$.
\end{tabular}
\right.
\end{equation*}
The numerical evidence of the optimality can be shown using the loading-unloading regime for the specimen. In terms of the force-displacement curve $\bar{F}$ vs.\ $\bar{u}_n$, where $\bar{F}$ is the reaction force at the top boundary, the first stage is the elastic loading, when no propagation occurs. It is then followed by the softening part, when the monotonically increasing applied displacement results in the crack advancement. The third stage is the elastic unloading. In this regime, provided the penalty parameter $\gamma$ is sufficiently large, one should be able to observe the linear decrease of $\bar{F}$ towards zero. More importantly, the appropriate $\gamma$ must provide a constant value for the regularized crack surface energy $\widetilde{E}_S$ in (\ref{Es_gam}).   

\begin{figure}[!ht]
\begin{center}
\includegraphics[width=1.0\textwidth]{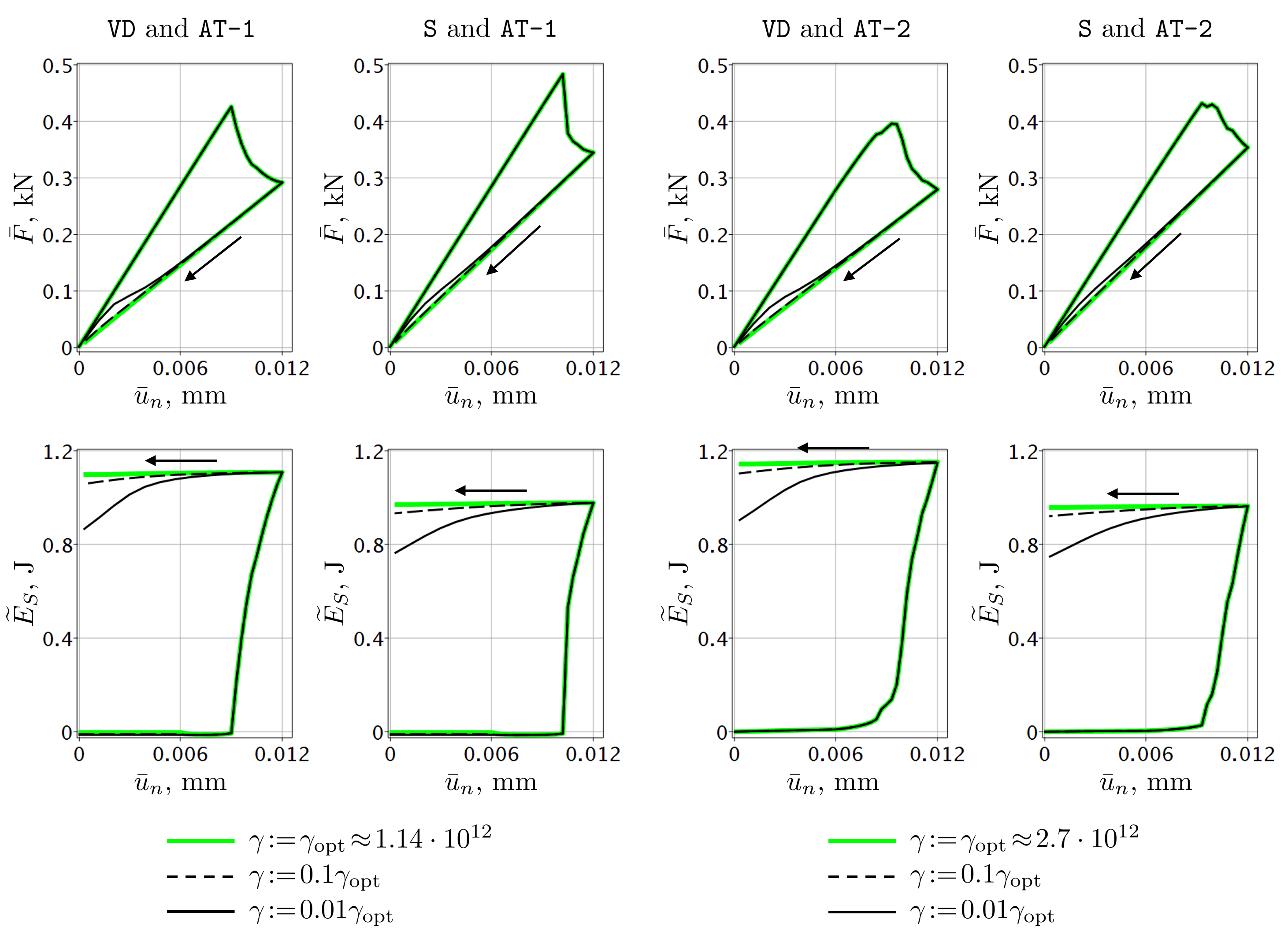}
\end{center}
\caption{Load-displacement and fracture surface energy-displacement curves for various combinations of split and model type, as well as for various magnitudes of the penalty parameter $\gamma$; the arrow designates the unloading branch.}
\label{Fig13}
\end{figure}

With $\bar{u}_1=6\cdot 10^{-3}$ mm as the displacement loading at the first step, and $\Delta\bar{u}:=0.3\cdot 10^{-3}$ mm as the loading increment, the adopted loading-unloading regime is defined by $\bar{u}_{n+1}:=\bar{u}_n+\Delta\bar{u}$ with $n=1,...,20$, and $\bar{u}_{n+1}:=\bar{u}_n-3\Delta\bar{u}$ with $n=21,...,33$. In Figure \ref{Fig13}, for every considered combination of split and model type, we plot the $\bar{F}$ vs.\ $\bar{u}_n$ and $\widetilde{E}_S$ vs.\ $\bar{u}_n$ curves. The arrow indicates the unloading branch. In each case, the dependence on the magnitude of $\gamma$ is evaluated. For the above optimal value $\gamma_\mathrm{opt}$, though derived using one-dimensional considerations, the expected correct trends for both $\bar{F}$ and $\widetilde{E}_S$ within the unloading phase are achieved. The first smaller value $0.1\gamma_\mathrm{opt}$ still provides the desired linear decrease of $\bar{F}$ similar to the optimal case, yet a decaying $\widetilde{E}_S$ is visible. For the second smaller value $0.01\gamma_\mathrm{opt}$, neither the behavior of $\bar{F}$, nor that of $\widetilde{E}_S$ are appropriate.

Figure \ref{Fig13} demonstrates that the intuitive choice of the penalty constant $\gamma$ is not feasible: all three values of $\gamma$ are large and of comparable orders of magnitude, but only the optimal one enforces the irreversibility with sufficient accuracy. It can also be concluded that, as claimed earlier, for a given fixed $\gamma_\mathrm{opt}$ the irreversibility behavior is not affected by the type of split. Finally, from Figure \ref{Fig13} the quantitative impact of tension-compression split and model type on the system response can be grasped. Let us note that  a detailed comparison of the formulations is out of scope here. The $\mathtt{AT}$-1 model, as expected, leads to a linear elastic limit before onset of fracture, whereas the $\mathtt{AT}$-2 one does not.

\subsubsection{Comparison of the irreversibility approaches in (\ref{Weak1}) and (\ref{Weak4})}
\begin{figure}[!ht]
\begin{center}
\includegraphics[width=1.0\textwidth]{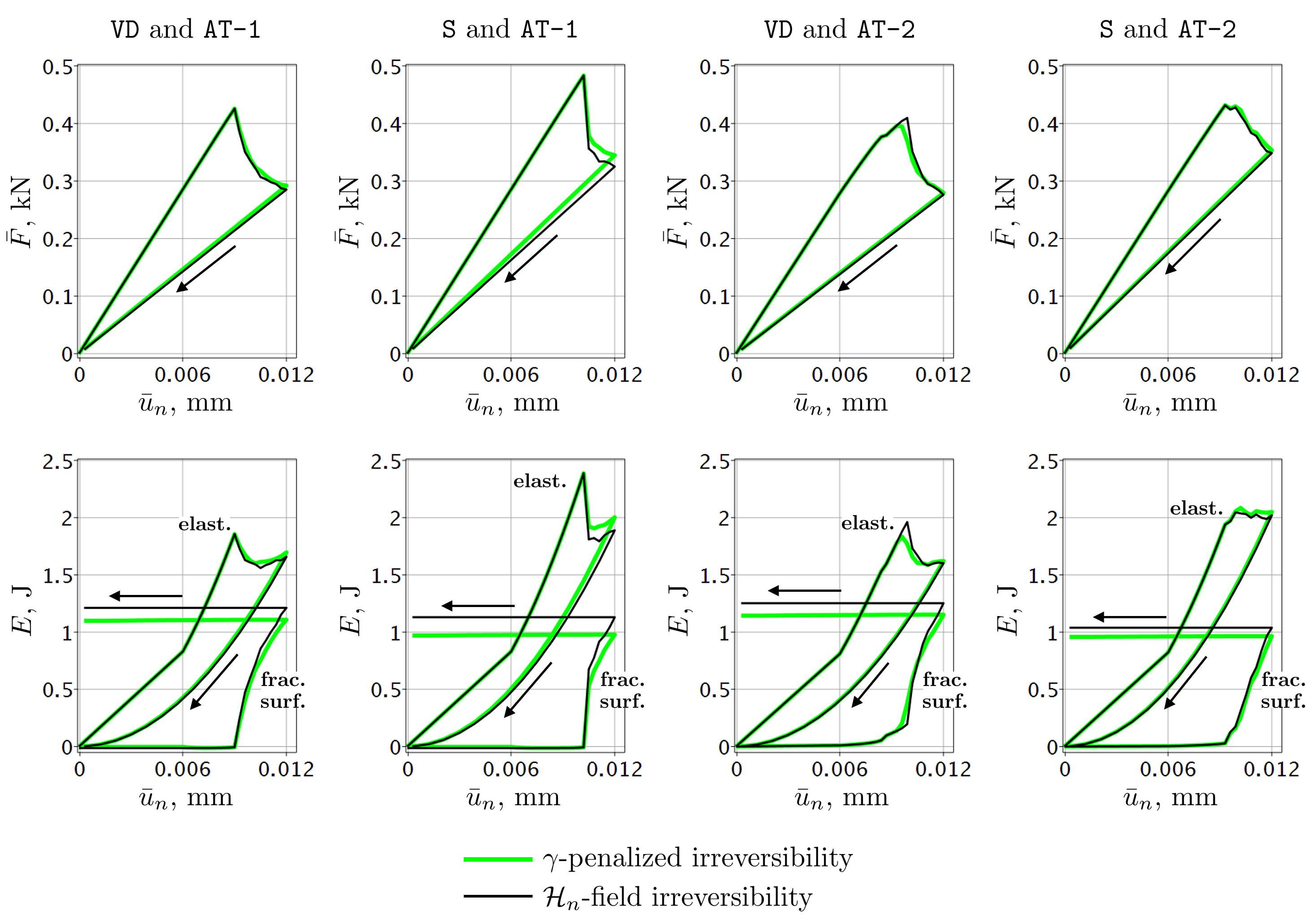}
\end{center}
\caption{Comparison of the load-displacement, elastic energy-displacement and fracture surface energy-displacement curves for the irreversibility approaches in (\ref{Weak1}) and (\ref{Weak4}).}
\label{Fig14}
\end{figure}

\begin{figure}[!ht]
\begin{center}
\includegraphics[width=1.0\textwidth]{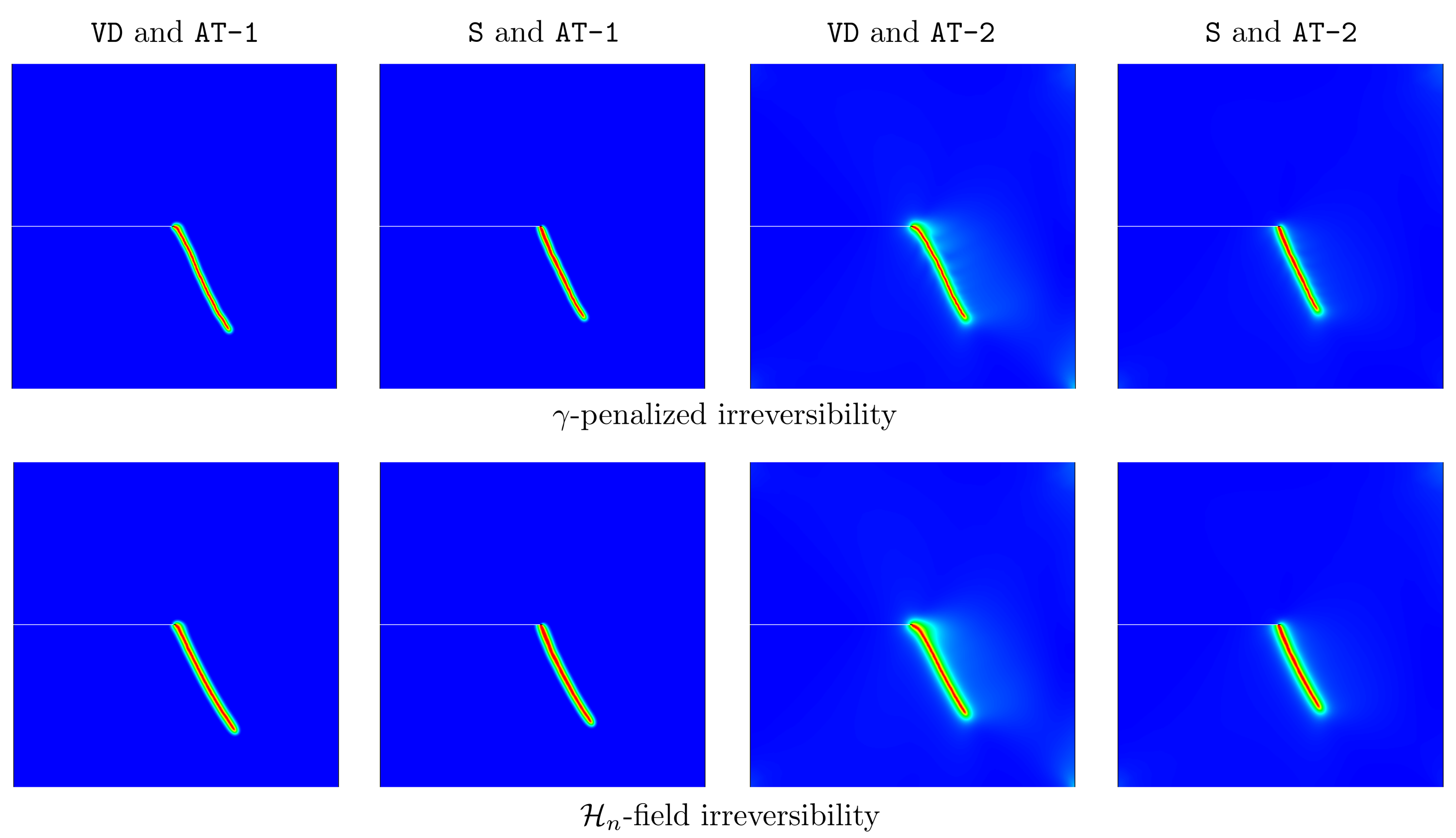}
\end{center}
\caption{Comparison of the crack phase-field profile at loading step $n=34$ for the irreversibility approaches in  (\ref{Weak1}) and (\ref{Weak4}).}
\label{Fig15}
\end{figure}

Herein, for the loading-unloading regime in section 4.1.1, we compare the penalized irreversibility results obtained for $\gamma_\mathrm{opt}$ with those obtained by using the history field technique originally introduced in Miehe et al.\ \cite{Miehe2010b}. In other words, for every fixed combination of split and model type we compare the solutions of the weak formulations (\ref{Weak1}) and (\ref{Weak4}) in terms of the load-displacement, elastic energy-displacement and fracture surface energy-displacement curves, as well as the crack phase-field profile at loading $\bar{u}_{34}=0.3\cdot 10^{-3}$ mm, see Figures \ref{Fig14} and \ref{Fig15}, respectively. Note that for computing the fracture surface energy in case of the $\gamma$-penalized formulation we use $\widetilde{E}_S$ as in (\ref{Es_gam}), and in the case of the history field formulation we use $E_S$ as in (\ref{Es}).

In all cases, both irreversibility approaches yield rather similar, but not identical results. Granting the $\gamma$-penalized results the reference, it can be seen that the history field approach yields under-estimation of the elastic energy and over-estimation of the fracture surface energy. We attribute this to the observation drawn from Figure \ref{Fig15} that the support of the phase-field profile in the latter case is visibly thicker. On the other hand, both approaches guarantee irreversibility, in the sense that during unloading they both lead to a linear decrease of the reaction force to zero as well as to a constant value of the crack surface energy. The natural conclusion is that, as argued in Section 2.3, formulation (\ref{Weak4}) is not equivalent to formulation (\ref{Weak1}), hence it is also not equivalent to the reference formulation in (\ref{Weak0}).

\subsection{Sneddon-Lowengrub problem}
\begin{figure}[!ht]
\begin{center}
\includegraphics[width=1.0\textwidth]{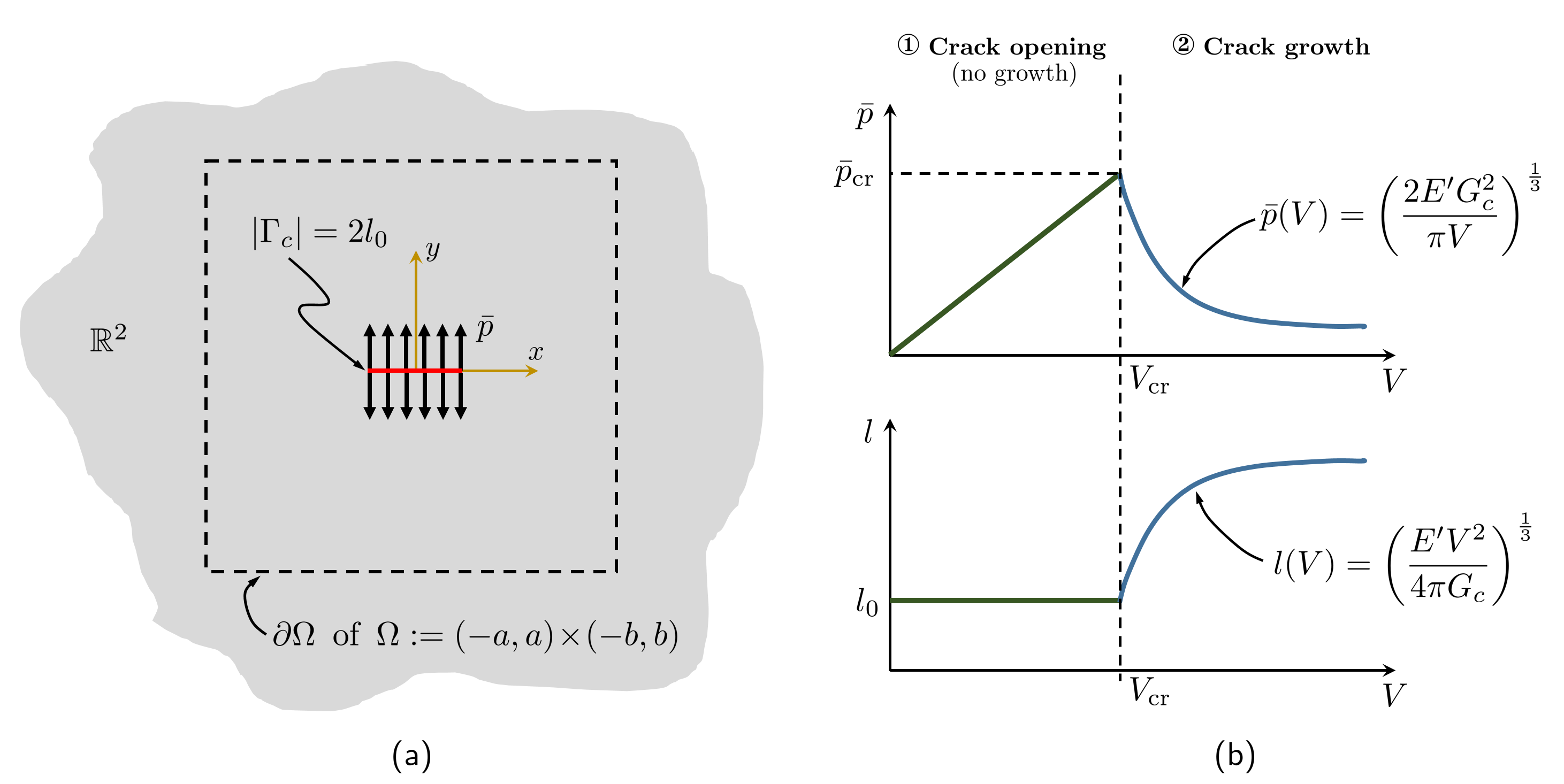}
\end{center}
\caption{(a) Geometry and loading setup for the Sneddon-Lowengrub problem; (b) evolution of the pressure $\bar{p}$ and fracture length $l$ as a function of the injected fluid volume $V$ depicted as a two-stage process.}
\label{Fig16}
\end{figure}

This problem considers an infinite domain $\mathbb{R}^2\backslash\Gamma_c$ with $\Gamma_c$ as a pre-existing crack of length $2l_0$ in the $y=0$ plane. Plane strain is assumed. Let $V$ be the volume of a fluid injected into the crack which generates the pressure $\bar{p}$ on the crack faces $\Gamma_c^\pm$. Also, let $V_\mathrm{cr}$ be the critical volume such that for all $V\leq V_\mathrm{cr}$ the crack opens but does not grow, and the pressure grows linearly until it reaches the critical value $\bar{p}_\mathrm{cr}$. When the injected volume exceeds $V_\mathrm{cr}$, the crack starts propagating and the pressure drops. The described situation is explained in detail in \cite{Bourdin2012hf} and is sketched in Figure \ref{Fig16}. Figure \ref{Fig16}(a) also depicts the computational domain $\Omega$ for the further fracture phase-field treatment of the problem. 
 
In the following, we are interested in the first stage of the crack evolution process in Figure \ref{Fig16}(b), that is, the crack opening stage. This case is typically referred to as the Sneddon-Lowengrub problem, since for any $\bar{p}\in(0,\bar{p}_\mathrm{cr}]$ the closed form solutions for the displacement and the stress fields in $\mathbb{R}^2\backslash\Gamma_c$ are obtained by Sneddon and Lowengrub in \cite{Sneddon1969}, see also \cite{Gudmundsson2014}. In particular, the $y$-component of the displacement $\bm u$ restricted to the upper crack face $\Gamma_c^+$ is given by
\begin{equation*}
u_y^+(x,0)=\frac{2\bar{p}l_0}{E^\prime}
\left(1-\frac{x^2}{l_0^2}\right)^\frac{1}{2}, \quad -l_0\leq x\leq l_0,
\end{equation*}
where $E^\prime:=\frac{E}{1-\nu^2}$, $E$ is the Young's modulus and $\nu$ is the Poisson's ratio of the material. Using the definition of the displacement jump across $\Gamma_c$, namely, $[\![\bm u]\!]_{\Gamma_c}:=u_y^+(x,0)-u_y^-(x,0)$, the anti-symmetry of $u_y$ with respect to the $x$-axis, and also $\bm n_{\Gamma_c}:=\bm n_{\Gamma_c^+}=(0,1)$, the COD reads
\begin{equation}
\mathrm{COD}(x):=[\![\bm u]\!]_{\Gamma_c}\cdot\bm n_{\Gamma_c}=2u_y^+(x,0).
\label{COD_Discr}
\end{equation}

For computational purposes we choose a rectangular domain $\Omega:=(-a,a)\times(-b,b)$ containing $\Gamma_c$. The phase-field formulation of the Sneddon-Lowengrub problem relies on the energy functional 
\begin{equation}
\mathcal{E}_\mathrm{HF}(\bm u,\alpha)=\mathcal{E}(\bm u,\alpha)
+\frac{\gamma}{2}\int_\Omega\langle \alpha-\alpha_{n-1}\rangle_{-}^2 \, \mathrm{d}{\bf x}
+\int_\Omega\bar{p}\bm u\cdot\nabla\alpha\, \mathrm{d}{\bf x},
\label{HF_PF}
\end{equation}
where the subscript HF stands for hydraulic fracture, $\mathcal{E}$ is given by (\ref{RegVAF1}) with the excluded traction contribution, and the second term is our penalty functional given by (\ref{OurPen}). The last term in (\ref{HF_PF}) is the energy term associated to the phase-field description of the pressure $\bar{p}\in(0,\bar{p}_\mathrm{cr}]$ applied to the crack faces, as originally proposed in \cite{Bourdin2012hf}. Given $\bar{p}$ and having found the solution $(\bm u,\alpha)\in {\bf V}_0\times H^1(\Omega)$ to the minimization problem for $\mathcal{E}_\mathrm{HF}$, the COD can be computed as \cite{Bourdin2012hf}:
\begin{equation}
\mathrm{COD}_\mathrm{PF}(x):=-\int_{-b}^b \bm u(x,y)
\!\cdot\!\nabla \alpha(x,y)\,\mathrm{d}y,
\label{COD_PF}
\end{equation}
This is an approximation to the COD in (\ref{COD_Discr}) which applies to the discrete crack setting. 

\begin{figure}[!ht]
\begin{center}
\includegraphics[width=1.0\textwidth]{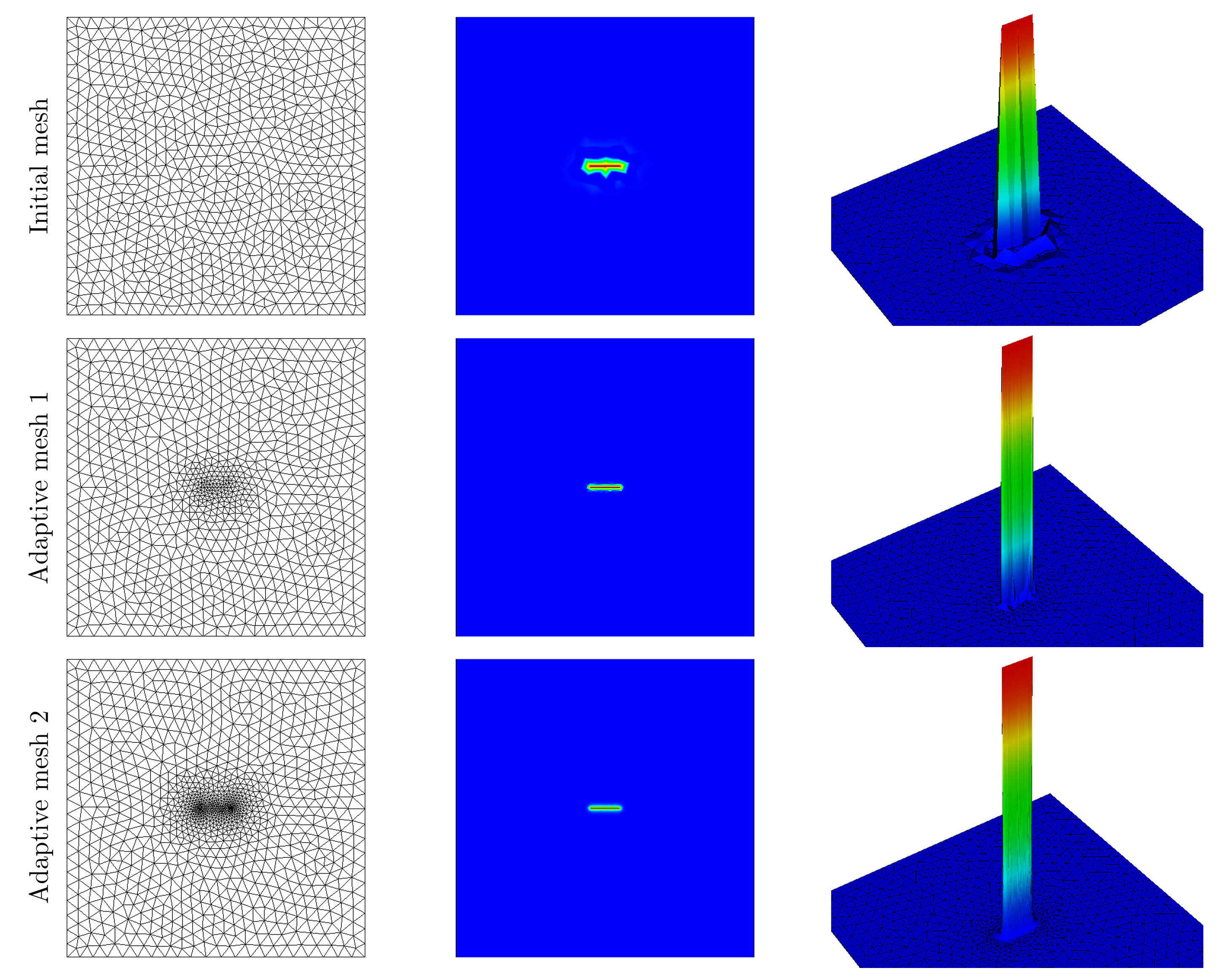}
\end{center}
\caption{Example of error-controlled adaptively refined meshes and the corresponding crack phase-field profile $\alpha_0$ computed on them.}
\label{Fig17}
\end{figure}

In our numerical experiment, the following dimensionless data are adopted from \cite{Bourdin2012hf}: $E=1$, $\nu=0.2$, $G_c=1$, $a=b=2$ and $l_0=0.2$. The length scale is set as follows $\ell:=\frac{2b}{200}=0.02$. Knowing also from \cite{Bourdin2012hf} that $\bar{p}_\mathrm{cr}:=\sqrt{G_cE^\prime/(\pi l_0) }\approx 1.288$, we take $\bar{p}=0.1$. Using (\ref{gamOpt}) with $\mathtt{TOL}_\mathrm{ir}:=0.01$, we obtain and set
\begin{equation*}
\gamma_\mathrm{opt}\approx
\left\{
\begin{tabular}{ll}
$2.1\cdot 10^{5}$, & when ${\sf w}(\alpha)=\alpha$,   \\
$5\cdot 10^{5}$,  & when ${\sf w}(\alpha)=\alpha^2$.
\end{tabular}
\right.
\end{equation*}
Finally, we want to compute the solution to (\ref{HF_PF}) in one loading step. In this case, in (\ref{HF_PF}), $\alpha_{n-1}:=\alpha_0$, where $\alpha_0$ is to be recovered by solving the minimization problem for
\begin{equation}
E_S(\alpha)=\left\{
\begin{tabular}{ll}
$\displaystyle 
\frac{3}{8}G_c \int_\Omega 
\left(\frac{\alpha}{\ell}+\ell|\nabla\alpha|^2\right)\mathrm{d}{\bf x}
+\frac{\rho}{2}\int_\Omega\langle \alpha \rangle_{-}^2 \; \mathrm{d}{\bf x}$, & when ${\sf w}(\alpha)=\alpha$, \\ [0.2cm]
$\displaystyle 
\frac{1}{2}G_c \int_\Omega 
\left(\frac{\alpha^2}{\ell}+\ell|\nabla\alpha|^2\right)\mathrm{d}{\bf x}$, & when ${\sf w}(\alpha)=\alpha^2$,
\end{tabular}
\right.
\label{Concr_Es}
\end{equation}
for all $\alpha\in H^1(\Omega)$ such that $\alpha(\Gamma_c)=1$. Note that this corresponds to modeling of the initial crack $\Gamma_c$ by the phase-field. To control the accuracy of $\alpha_0$, an error-controlled adaptive mesh refinement strategy developed in our forthcoming paper \cite{Gerasimov2018b} is used. As an example, Figure \ref{Fig17} depicts the sequence of adaptive meshes and the corresponding $\alpha_0$ computed on them for the case of the quadratic model. The total amount of adaptive steps we perform is five, and the figure presents only two of them. The final adaptive mesh is such that the minimal element size in the localization zone of $\alpha$ is $h\approx 7.6\cdot 10^{-5}$ (for the linear model, it is $h\approx 7.3\cdot 10^{-5}$), that is, $h\ll\ell$ is fulfilled.

\begin{figure}[!ht]
\begin{center}
\includegraphics[width=1.0\textwidth]{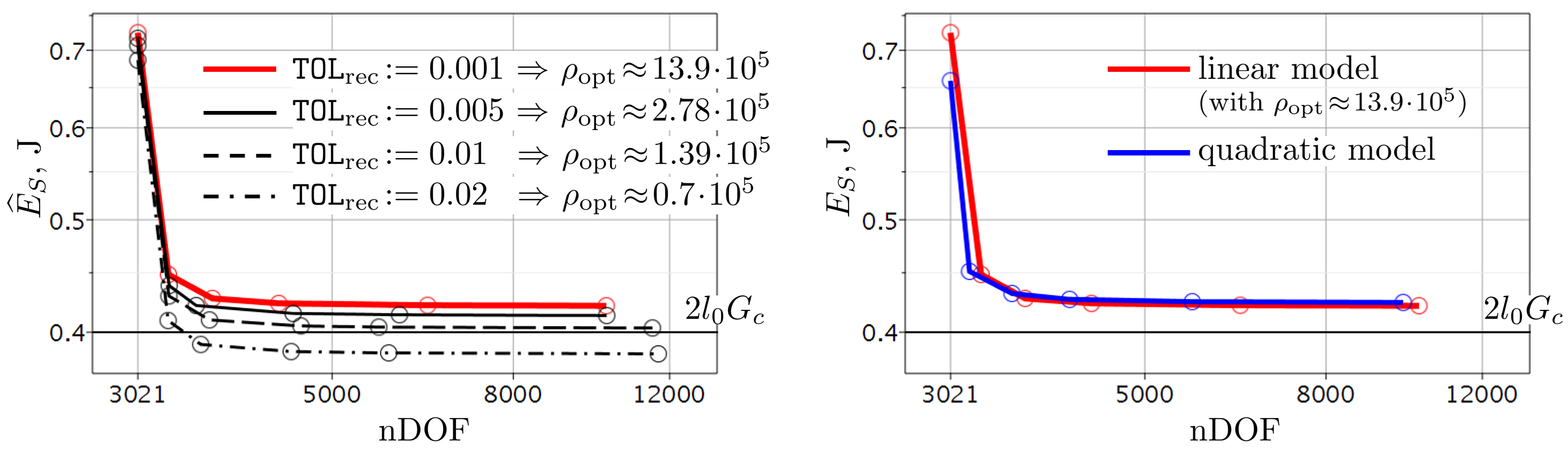}
\end{center}
\caption{On the left: convergence of the fracture surface energy $E_S$ for the linear model on the adaptive meshes for various magnitudes of the recovery tolerance threshold $\mathtt{TOL}_\mathrm{rec}$ and the induced penalty parameter $\rho_\mathrm{opt}$; on the right: comparison of the two models.}
\label{Fig18}
\end{figure}

\begin{figure}[!ht]
\begin{center}
\includegraphics[width=1.0\textwidth]{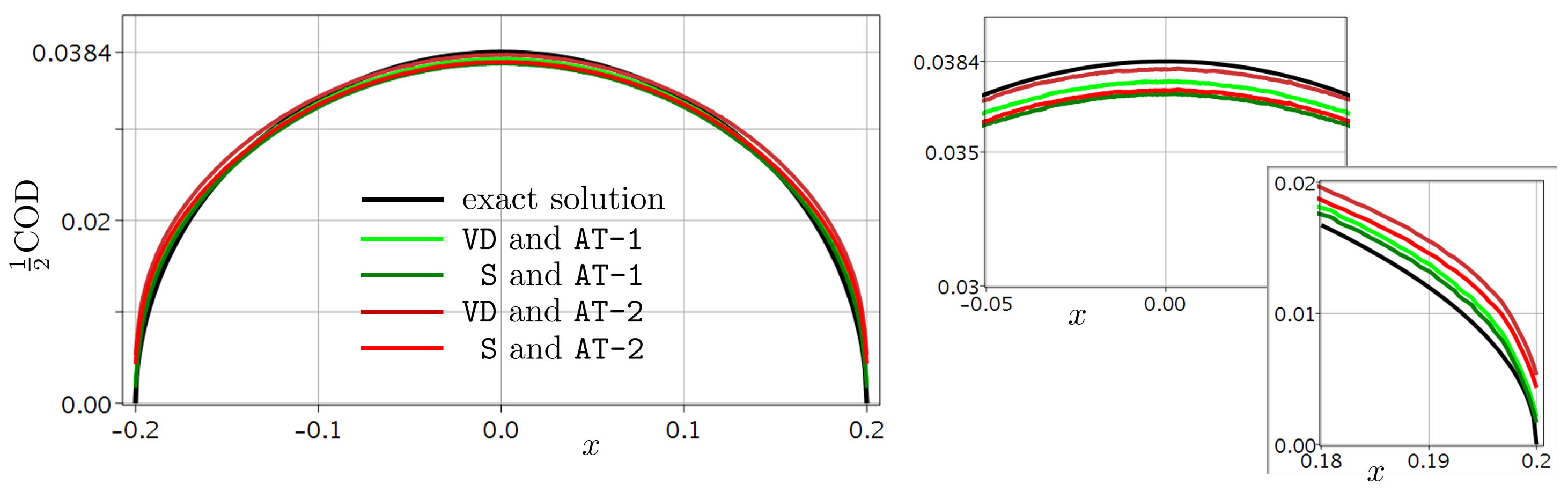}
\end{center}
\caption{Comparison of the COD obtained by (\ref{COD_PF}) for various combinations of split and model type with the exact Sneddon-Lowengrub solution in (\ref{COD_Discr}).}
\label{Fig19}
\end{figure}

Using the adaptively refined meshes, two kinds of convergence behavior of the fracture surface energy $E_S$ for the linear model can be observed: the finite element convergence and the convergence related to the penalty parameter $\rho$, see the left plot of Figure \ref{Fig18}. For any fixed magnitude of the recovery tolerance threshold $\mathtt{TOL}_\mathrm{rec}$ and the corresponding $\rho_\mathrm{opt}$, an asymptotic convergence of $E_S$ from above in terms of the degrees of freedom to a horizontal limit is visible. It can also be seen that the asymptotic limit of $E_S$ increases with the increase of $\rho_\mathrm{opt}$, as expected, and has a supremum when $\rho\rightarrow+\infty$. With the given finite length scale $\ell$, this supremum overestimates the reference $\Gamma$-convergence value $G_c|\Gamma_c|=2l_0G_c$ (this also holds for the quadratic model, see the right plot of Figure \ref{Fig18})\footnote{Such an overestimation occurs in the first place due to the finiteness of $\ell$. Secondly, the relative size of the half-spherical parts of the support of $\alpha$ around the crack tips -- in the sketch below these are denoted as $T_1$ and $T_2$ -- with respect to the size of $\Omega_{\Gamma_c}$ contributes as well. More precisely, for fixed $\ell$, the {\em lengthier} the part $\Omega_{\Gamma_c}$ is, the {\em less} pronounced overestimation will occur. This trend represented in terms of the crack length $|\Gamma_c|$ is reported in the corresponding table. Therein, $E_S$ is the converged value after 5 adaptive re-meshing steps, $\ell=0.02$ and $l_0=0.2$ are our reference values, for which the results presented in Figures \ref{Fig17}, \ref{Fig18} and \ref{Fig19} are computed.

\begin{minipage}{0.3\textwidth}
	\begin{center}
	\resizebox{3.0cm}{!} {\includegraphics{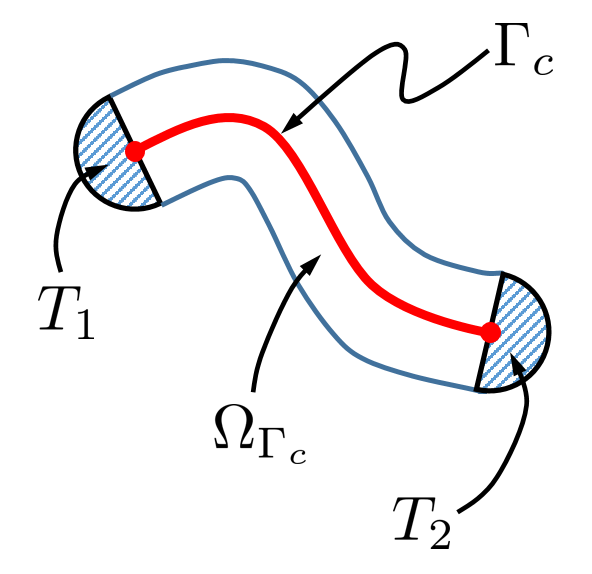}}
	\end{center}
\end{minipage}
\begin{minipage}{0.7\textwidth}
	\begin{center}
	\begin{tabular}{c|c|c}
      \hline
       Model &  $|\Gamma_c|$   &  $\frac{E_S-G_c|\Gamma_c|}{G_c|\Gamma_c|}\cdot 100\%$  \\
	\hline
       \begin{tabular}{c} linear \\  (with $\rho_\mathrm{opt}\approx13.9\cdot10^5$) \end{tabular}  
       		&  \begin{tabular}{c} $2l_0$ \\  $4l_0$ \\ $8l_0$ \end{tabular}  
			&  \begin{tabular}{c} $5.557$ \\ $3.205$ \\ $1.883$ \end{tabular} \\
	\hline
	quadratic &  \begin{tabular}{c} $2l_0$ \\  $4l_0$ \\ $8l_0$ \end{tabular}  
			 &  \begin{tabular}{c} $6.305$ \\ $3.741$ \\ $2.108$ \end{tabular} \\
	\hline
	\end{tabular}
	\end{center}
\end{minipage}
It can also be observed that, as expected, the linear model yields a better crack approximation than the quadratic one.}. In contrast, insufficiently large values of $\rho$ cannot enforce the condition $\alpha\geq0$ accurately, thus resulting in underestimation of $2l_0G_c$, as indicated by the dash-and-dot curve in the left plot. It can clearly be concluded that the accurate intuitive choice of the penalty parameter $\rho$ is not feasible. 

Finally, in Figure \ref{Fig19} we present the COD computed by (\ref{COD_PF}) for various combinations of split and model type and compare these results with the exact Sneddon-Lowengrub solution given by (\ref{COD_Discr}). As can be observed in the main plot, for the chosen optimal penalty parameters $\gamma_\mathrm{opt}$, all four phase-field formulations approximate the exact solution accurately enough. They also yield rather similar results. To outline the differences, the corresponding regions -- the crack face midspan and the crack tip vicinity -- are zoomed in. In the first region, the crack phase-field solution underestimates the analytical one, whereas in the second region the situation is the opposite and such that $\mathrm{COD}_\mathrm{PF}(\pm l_0)\neq0$. Similar findings are presented e.g.\ in \cite[section 4.4]{Wilson2016}\footnote{The formulation in \cite{Wilson2016} uses $\mathtt{S}$-split, ${\sf g}(\alpha):=k\left(1-\left(\frac{k-1}{k} \right)^{(1-\alpha)^2}\right)$, $k>0$, ${\sf w}(\alpha):=\alpha^2$ and a more complicated representation for the pressure term than the one in (\ref{HF_PF}).}. This can be attributed to the difficulty to obtain both solution fields accurately within the finite element setting. It can particularly be noticed that the $\mathtt{AT}$-1 model (regardless of split type) is superior to the $\mathtt{AT}$-2 formulation in the vicinity of the crack tip. Indeed, no premature evolution of $\alpha$ in the elastic loading regime occurs when the $\mathtt{AT}$-1 model is used. 

\section{Conclusions}

In the phase-field model of brittle fracture, irreversibility of the crack phase-field imposed to prevent fracture healing leads to a constrained minimization problem. In this paper, we focused on treatment of the irreversibility constraint via a simple penalty formulation, which, provided the penalty constant is well-tuned, is a good approximation to the original one. Exploiting the notion of the optimal phase-field profile as well as the $\Gamma$-convergence result, we proposed an analytical procedure in the one-dimensional setting for deriving a lower bound for the penalty constant which guarantees a sufficiently accurate enforcement of the crack phase-field irreversibility. This lower bound was found to be a function of two formulation parameters (the fracture toughness and the regularization length scale) and to be independent on the problem setup (geometry, boundary conditions etc.) and on the formulation ingredients (degradation function and tension-compression split). For the recovery of the phase-field profile within a formulation using a linear damage function, we similarly addressed penalization as a means of enforcing the non-negativity of the phase-field variable. The optimally-penalized formulation for both irreversibility and phase-field recovery was tested for two benchmark problems in two dimensions, including one with available analytical solution. The numerical results demonstrated that the penalty parameter computed with the developed relationship was the smallest value able to enforce the constraint with sufficient accuracy. We also compared our results for irreversibility with those obtained by using the notion of the history field by Miehe et al.\ (2010). The comparison showed that the approach based on the history field is able to enforce irreversibility, however it leads to the solution of a problem that is not equivalent to the original one.

Due to its generality, the presented approach can be extended to different phase-field models of fracture. E.g. it can be readily applied to derive a lower bound for the penalty parameter in the penalized formulation that uses the double-well ${\sf w}(\alpha)=16\alpha^2(1-\alpha)^2$ and the double-obstacle ${\sf w}(\alpha)=4\alpha(1-\alpha)$ functions. These are typical ingredients of fracture phase-field models stemming from the physics community, see e.g.\ \cite{Kuhn2015} and the review paper \cite{Ambati2015review}. In terms of the optimal phase-field profile, the double-well function is similar to the quadratic model, wheres the double-obstacle to the linear one. These and other extensions may well be considered as future developments of this research.

{\small
\paragraph{Acknowledgment.} The support of the DFG through the project `Phase-field computation of brittle fracture: robustness, efficiency, and characterisation of solution non-uniqueness' is gratefully acknowledged. Some of the presented results were obtained and the part of the manuscript was written during the 5 month research stay of the first author in the group of Prof. Andrew Chan, the Colledge of Science and Engineering, the University of Tasmania, Australia. This visit was funded by an Australia Award -- 2018 Endeavour Research Fellowship.} 

\appendix
\section{Solution to problem (\ref{bvp1d_Lin})}
\label{Ap0}
\begin{figure}[!ht]
\begin{center}
\includegraphics[width=1.0\textwidth]{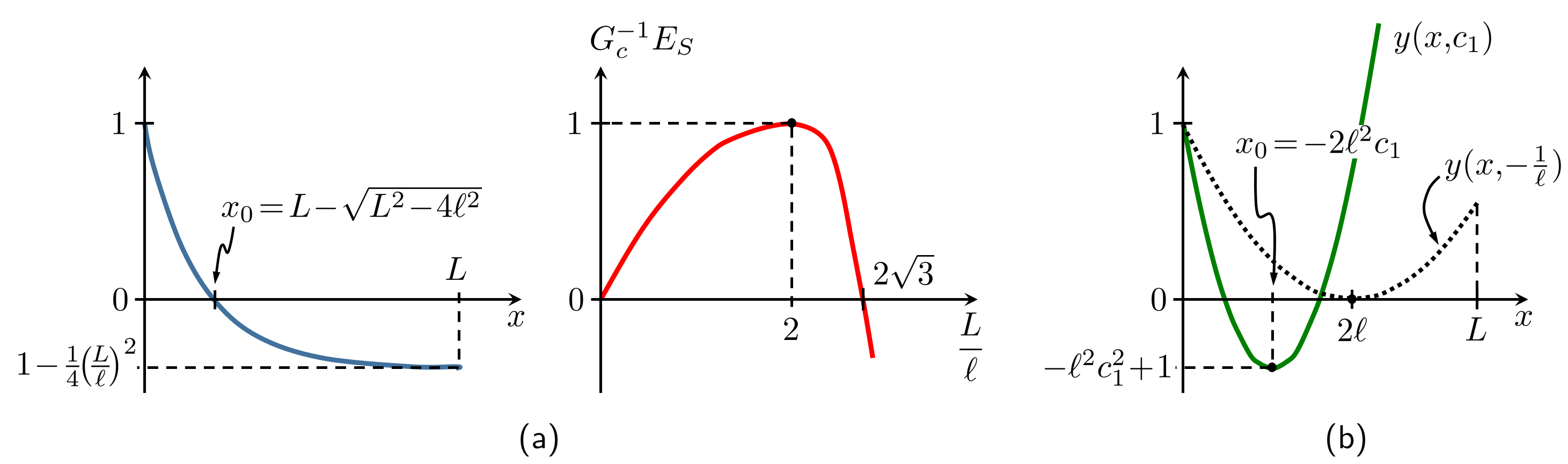}
\end{center}
\caption{(a) On the left: minimizer $\alpha^\star$ of the linear model without using the constraint $\alpha\geq0$, it is negative in $(x_0,L)$ if $\frac{L}{\ell}>2$; on the right: the plot of $G_c^{-1}E_S(\alpha^\star)$ as a function of $\frac{L}{\ell}$, the reference value $1$ is attained when $\frac{L}{\ell}=2$; (b) Stages for constructing the minimizer $\alpha^\star$ given by (\ref{aStar1d_Lin}).}
\label{Fig20}
\end{figure}

Let us first show that without the constraint $\alpha\geq0$, the correct optimal profile $\alpha^\star$ cannot be constructed in this case. Indeed, omitting $\alpha\geq0$, one arrives at the boundary value problem (\ref{bvp1d_Lin}) with the inequalities replaced by equalities. In this situation, the general solution of the differential equation reads
\begin{equation*}
\alpha(x)=\frac{1}{4}\frac{x^2}{\ell^2}+c_1x+c_2, \quad c_1,c_2\in\mathbb{R}.
\end{equation*}
Using the boundary conditions, one obtains $c_1=-\frac{1}{2}\frac{L}{\ell^2}$ and $c_2=1$, yielding
\begin{equation*}
\alpha^\star(x)=\frac{1}{4}\frac{x^2}{\ell^2}-\frac{1}{2}\frac{L}{\ell^2}x+1,
\end{equation*}
that combined with (\ref{Es1d}) yields
\begin{equation*}
G_c^{-1}E_S(\alpha^\star)=-\frac{1}{16}\left(\frac{L}{\ell}\right)^3+\frac{3}{4}\frac{L}{\ell}.
\end{equation*}
It can immediately be noticed that $\alpha^\star(L)=1-\frac{1}{4}\left(\frac{L}{\ell}\right)^2$ is negative when $\ell<\frac{L}{2}$, which is plausible owing to our main assumption $0<\ell\ll L$. We depict this in the left plot of Figure \ref{Fig20}(a). Therein, the point $x_0=L-\sqrt{L^2-4\ell^2}$ is the solution of $\alpha^\star(x)=0$. We also depict $G_c^{-1}E_S(\alpha^\star)$ as a function of $\frac{L}{\ell}$ in the right plot of Figure \ref{Fig20}(a). The desired value $1$ is attained when $\frac{L}{\ell}=2$, and $G_c^{-1}E_S$ becomes negative when $\frac{L}{\ell}>2\sqrt{3}$. Thus we need the introduction of the constraint $\alpha\geq0$ into the corresponding admissible set for the linear model $E_S$ in (\ref{Es1d}). 

We now address the actual formulation (\ref{bvp1d_Lin}). The solution of the differential inequality in (\ref{bvp1d_Lin}) is given by
\begin{equation*}
0\leq \alpha(x) \leq
\frac{1}{4}\frac{x^2}{\ell^2}+c_1x+c_2, \quad c_1,c_2\in\mathbb{R}.
\end{equation*}
Using the boundary condition $\alpha(0)=1$, one obtains $c_2=1$ and arrives at the condition which the candidate for $\alpha^\star$ must fulfill: 
\begin{equation}
0\leq \alpha(x) \leq \frac{1}{4}\frac{x^2}{\ell^2}+c_1x+1 \quad \mathrm{in} \; (0,L),
\label{sys}
\end{equation}
For obtaining the constant $c_1$, we define $\frac{1}{4}\frac{x^2}{\ell^2}+c_1x+1=:y(x,c_1)$. It is a parabola with the vertex of abscissa $x_0:=-2\ell^2c_1$ and such that $y(x_0,c_1)=-\ell^2c_1^2+1$, see Figure \ref{Fig20}(b). The demand $x_0\in(0,L)$ implies that $c_1<0$. Furthermore, the requirement $0\leq\alpha(x)\leq y(x,c_1)$ for all $x\in(0,L)$ implies in particular that $y(x_0,c_1)\geq0$, that is, $-\ell^2c_1^2+1\geq0$. In the limiting case, we may simply impose $-\ell^2c_1^2+1=0$, that yields $c_1=\pm\frac{1}{\ell}$. We take the negative value, thus ultimately setting $c_1:=-\frac{1}{\ell}$. This, in turn, yields $x_0=2\ell$ and $y(x,-\frac{1}{\ell})=\frac{1}{4}\frac{x^2}{\ell^2}-\frac{1}{\ell}x+1$, see the dashed line in Figure \ref{Fig20}(b). As a result, system (\ref{sys}), defining the admissible set of minimizers, turns into
\begin{equation}
0\leq \alpha(x) \leq \frac{1}{4}\frac{x^2}{\ell^2}-\frac{1}{\ell}x+1 \quad \mathrm{in} \; (0,L).
\label{sys1}
\end{equation}
Clearly, any $\alpha$ satisfying (\ref{sys1}) also satisfies the boundary condition $\alpha^\prime(L)\geq0$ in (\ref{bvp1d_Lin}), as required. 

The most straightforward way to choose only one representative from (\ref{sys1}) is to set $\alpha^\star$ to the function $y(x,-\frac{1}{\ell})$ in the interval $(0,2\ell)$, whereas in the interval $[2\ell,L)$, $\alpha^\star$ must be set to $0$ in order to comply with the crack phase-field formalism. As a result, we obtain
\begin{equation*}
\alpha^\star(x)=
\left\{
\begin{tabular}{cl}
$\displaystyle \frac{1}{4}\frac{x^2}{\ell^2}-\frac{x}{\ell}+1
=\left(1- \frac{x}{2\ell}\right)^2$, & in $(0,2\ell)$.   \\[0.2cm]
$0$, & in $[2\ell,L)$.
\end{tabular}
\right.
\end{equation*}

\section{Solution to problem (\ref{bvp_s})}
\label{Ap1}
In either case of linear and quadratic model, problem (\ref{bvp_s}) is a semi-linear problem and requires an iterative solution procedure. We use the standard incremental setting: at every iteration $i\geq0$, given known $\alpha_i$ and inserting $\alpha_{i+1}:=\alpha_i+\Delta\alpha_i$ into (\ref{bvp_s}), the linearized problem for the unknown $\Delta\alpha_i$ is obtained:
\begin{equation}
\left\{
\begin{tabular}{rl}
$\displaystyle -\ell^2(\Delta\alpha_i)^{\prime\prime}
+\left[ \frac{1}{2}{\sf w}^{\prime\prime}(\alpha_i)+s{\sf H}^{-}(\alpha_i-\alpha^\star) \right] \Delta\alpha_i
=\ell^2\alpha_i^{\prime\prime}-\frac{1}{2}{\sf w}^{\prime}(\alpha_i)-s \langle \alpha_i-\alpha^\star \rangle_{-}$ & in $(0,L)$,   \\ 
$(\Delta\alpha_i)^\prime=-\alpha_i^\prime$ & at $x\in\{0,L\}$.
\end{tabular}
\right.
\label{bvplin_s}
\end{equation}
Note that in the above we used $\langle \alpha_{i+1}-\alpha^\star \rangle_{-}=
\langle \alpha_i+\Delta\alpha_i-\alpha^\star \rangle_{-}
\approx \langle \alpha_i-\alpha^\star \rangle_{-}+
{\sf H}^{-}(\alpha_i-\alpha^\star)\Delta\alpha_i$, with ${\sf H}^{-}$ given by (\ref{Hmin}), and also 
${\sf w}^\prime(\alpha_{i+1})={\sf w}^\prime(\alpha_i+\Delta\alpha_i)
\approx {\sf w}^\prime(\alpha_i)+{\sf w}^{\prime\prime}(\alpha_i)\Delta\alpha_i$. Once $\Delta\alpha_i$ is known, the residual to be checked for the updated solution $\alpha_{i+1}$ at $i\geq0$ reads: 
\begin{equation}
\mathrm{Res}_i:=-\ell^2\alpha_{i+1}^{\prime\prime}+\frac{1}{2}{\sf w}^{\prime}(\alpha_{i+1})
+s \langle \alpha_{i+1}-\alpha^\star \rangle_{-}.
\label{res_s}
\end{equation}
Further developments for the two model types are presented below.

\subsection{Linear model}
\label{Ap1_LinMod}
The exact solution to problem (\ref{bvp_s}) when ${\sf w}(\alpha)=\alpha$ can be constructed using the linearized formulation (\ref{bvplin_s}) in only two iterations $i\geq0$. More precisely, already at iteration $i=1$ we obtain $\alpha_2=\alpha_1+\Delta\alpha_1=\alpha_0+\Delta\alpha_0+\Delta\alpha_1$, where the initial guess $\alpha_0$ is taken as $\alpha^\star$ in (\ref{aStar1d_Lin}), and the computed $\Delta\alpha_0$ and $\Delta\alpha_1$ are given by
\begin{equation}
\Delta\alpha_0(x)=\frac{1}{\sqrt{s}}\frac{1}{\exp\left( -2\sqrt{s}\frac{L}{\ell} \right)-1}
\left[
\exp\left( -\sqrt{s}\frac{x}{\ell} \right)
+\exp\left( -2\sqrt{s}\frac{L}{\ell} \right) \exp\left( \sqrt{s}\frac{x}{\ell} \right)
\right],
\label{Da0_Lin}
\end{equation}
and
\begin{equation}
\Delta\alpha_1(x)=\left\{
\begin{tabular}{l}
$\displaystyle \frac{1}{4s}\frac{\exp\left( -2\sqrt{s} \right)-\exp\left( -2\sqrt{s}\left(\frac{L}{\ell}-1\right) \right)}
							{\exp\left( -2\sqrt{s}\frac{L}{\ell} \right)-1}
			\left[
			\exp\left( -\sqrt{s}\frac{x}{\ell} \right)+\exp\left( \sqrt{s}\frac{x}{\ell} \right)
			\right]$,   \\ 
\hspace{9.0cm} in $(0,2\ell)$, \\ [0.2cm]
$\displaystyle \frac{1}{4s}\frac{\exp\left( -2\sqrt{s} \right)-\exp\left( 2\sqrt{s} \right)}
							{\exp\left( -2\sqrt{s}\frac{L}{\ell} \right)-1}
			\left[
			\exp\left( -\sqrt{s}\frac{x}{\ell} \right)
			+\exp\left( -2\sqrt{s}\frac{L}{\ell} \right) \exp\left( \sqrt{s}\frac{x}{\ell} \right)
			\right]-\frac{1}{2s}$,  \\
\hspace{10.7cm} in $[2\ell,L)$.
\end{tabular}
\right.
\label{Da1_Lin}
\end{equation}
Note that the piecewise function $\Delta\alpha_1$, similarly to the corresponding $\alpha^\star$,  is continuously-differentiable at $x=2\ell$. It is then straightforward to see that $\mathrm{Res}_1\equiv0$ indeed holds. The resulting solution $\alpha^\star_\gamma$ then reads:
\begin{equation*}
\alpha^\star_\gamma(x)=\alpha^\star(x)+\Delta\alpha_0(x)+\Delta\alpha_1(x).
\end{equation*}
The plots of the constructed $\alpha^\star_\gamma$ are depicted in Figure \ref{Fig4} in section \ref{LinMod_gam}. Using the assumption that $s\gg1$ and $\frac{L}{\ell}\geq2$, thus yielding $\exp\left(-\sqrt{s}\frac{L}{\ell} \right)\ll1$, the simplified versions of (\ref{Da0_Lin}) and (\ref{Da1_Lin}) can be obtained:
\begin{equation*}
\Delta\alpha_0(x)\approx -\frac{1}{\sqrt{s}} \exp\left( -\sqrt{s}\frac{x}{\ell} \right),
\end{equation*}
and
\begin{equation*}
\Delta\alpha_1(x)\approx \left\{
\begin{tabular}{ll}
$\displaystyle -\frac{1}{4s} \exp\left( -2\sqrt{s} \right)
			\left[
			\exp\left( -\sqrt{s}\frac{x}{\ell} \right)+\exp\left( \sqrt{s}\frac{x}{\ell} \right)
			\right]$,   &  in $(0,2\ell)$, \\ [0.5cm]
$\displaystyle -\frac{1}{4s} \left( \exp\left( -2\sqrt{s} \right)-\exp\left( 2\sqrt{s} \right) \right)
			\exp\left( -\sqrt{s}\frac{x}{\ell} \right)-\frac{1}{2s}$,  & in $[2\ell,L)$.
\end{tabular}
\right.
\end{equation*}
This proves the $\frac{L}{\ell}$-independence of $\alpha^\star_\gamma$ observed in Figure \ref{Fig4}. Interestingly, the simplified $\Delta\alpha_1$ remains continuously differentiable at $x=2\ell$ as its original counterpart.

\subsection{Quadratic model}
\label{Ap1_SqMod}
The exact solution to problem (\ref{bvp_s}) with ${\sf w}(\alpha)=\alpha^2$ can be constructed using the linearized formulation (\ref{bvplin_s}) in only one iteration. Thus, already at iteration $i=0$, we obtain $\alpha_1=\alpha_0+\Delta\alpha_0$, where we take the exact $\alpha^\star$ in (\ref{aStar1d_Sq}) as the initial guess $\alpha_0$. The computed $\Delta\alpha_0$ reads
\begin{align}
\Delta\alpha_0(x)= & \frac{1}{\sqrt{s+1}} \frac{1}{\exp\left( {-2\sqrt{s+1}\frac{L}{\ell}} \right)-1}
			    \frac{1-\exp\left( {-2\frac{L}{\ell}} \right)}{1+\exp\left( {-2\frac{L}{\ell}} \right)} \nonumber \\
\times &
\left[
\exp\left( -\sqrt{s+1}\frac{x}{\ell} \right)
+\exp\left( -2\sqrt{s+1}\frac{L}{\ell} \right) \exp\left( \sqrt{s+1}\frac{x}{\ell} \right)
\right].
\label{Da0_Sq}
\end{align}
It is straightforward to check that $\mathrm{Res}_0\equiv0$ holds. The resulting representation for $\alpha^\star_\gamma$ then follows:
\begin{equation*}
\alpha^\star_\gamma(x)=\alpha^\star(x)+\Delta\alpha_0(x).
\end{equation*}
The plot of the obtained $\alpha^\star_\gamma$ is presented in Figure \ref{Fig7} in section \ref{SqMod_gam}. Using again that $s\gg1$ and $\frac{L}{\ell}\geq2$, it can shown that the dependence of $\alpha^\star_\gamma$ on $\frac{L}{\ell}$ is negligible. Indeed, $\exp\left(-\sqrt{s+1}\frac{L}{\ell}\right)\ll1$ and $\exp\left(-2\frac{L}{\ell}\right)\ll1$ such that the simplified representation for (\ref{Da0_Sq}) is as follows:
\begin{equation*}
\Delta\alpha_0(x) \approx -\frac{1}{\sqrt{s+1}} \exp\left( -\sqrt{s+1}\frac{x}{\ell} \right).
\end{equation*}
This confirms the independence of $\alpha^\star_\gamma$ on $\frac{L}{\ell}$ observed in Figure \ref{Fig7}.

\section{Solution to problem (\ref{bvp_r})}
\label{Ap2}
\begin{figure}[!ht]
\begin{center}
\includegraphics[width=1.0\textwidth]{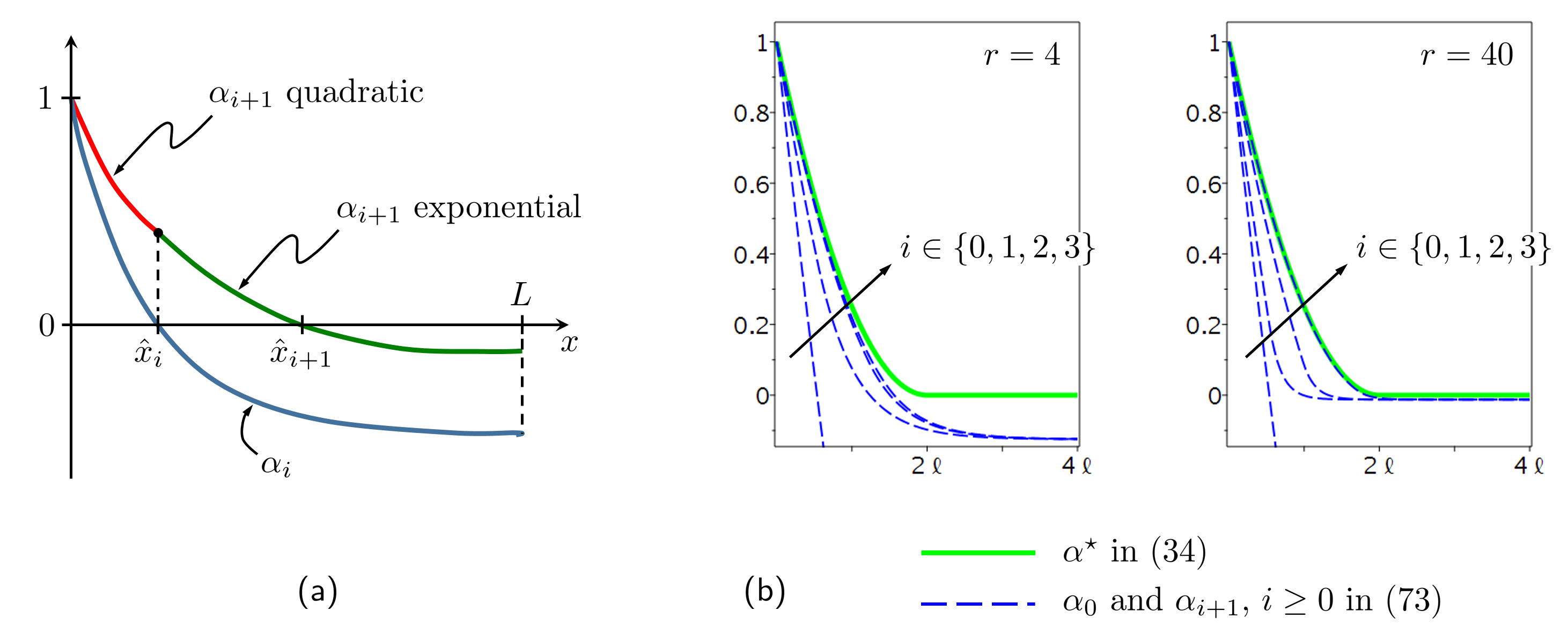}
\end{center}
\caption{(a) The construction idea of $\alpha_{i+1}$ in (\ref{ai_r}) using $\alpha_i$ and the corresponding point $\hat{x}_i$; (b) Convergence behavior of the sequence $\{\alpha_{i+1}\}$, $i\geq0$ in (\ref{ai_r}) for two different values of the penalty parameter $r$.}
\label{Fig21}
\end{figure}

The solution of the semi-linear problem in (\ref{bvp_r}) is constructed iteratively. The linearized equation to be solved for $\Delta\alpha_i$ at any iteration $i\geq0$ reads in this case
\begin{equation}
-\ell^2(\Delta\alpha_i)^{\prime\prime}
+r{\sf H}^{-}(\alpha_i) \Delta\alpha_i
=\ell^2\alpha_i^{\prime\prime}-\frac{1}{2}-r \langle \alpha_i \rangle_{-}
\quad \mathrm{in} \; (0,L),
\label{eqlin_r}
\end{equation}
It can be seen that for any $\alpha_i$, $i\geq0$ satisfying $\alpha_i(0)=1$, $\alpha_i^\prime(L)=0$ such that $\hat{x}_i\in(0,L)$ is a solution to $\alpha_i(x)=0$ -- see Figure \ref{Fig21}(a) for illustration -- equation (\ref{eqlin_r}) turns into
\begin{equation*}
\left\{
\begin{tabular}{ll}
$\displaystyle -\ell^2(\Delta\alpha_i)^{\prime\prime}=
\ell^2(\alpha_i)^{\prime\prime}-\frac{1}{2}$ & in $(0,\hat{x}_i)$,   \\ [0.2cm]
$\displaystyle -\ell^2(\Delta\alpha_i)^{\prime\prime}+r\Delta\alpha_i=
\ell^2(\alpha_i)^{\prime\prime}-\frac{1}{2}-r\alpha_i$ & in $(\hat{x}_i,L)$,
\end{tabular}
\right.
\end{equation*}
whose solutions are as follows:
\begin{equation*}
\Delta\alpha_i(x):=\left\{
\begin{tabular}{ll}
$\displaystyle -\alpha_i(x)+\frac{1}{4\ell^2}x^2+c_1x+c_2$, & in $(0,\hat{x}_i)$,   \\ [0.2cm]
$\displaystyle -\alpha_i(x)+b_1\exp\left(-\sqrt{r}\frac{x}{\ell}\right) + b_2\exp\left(\sqrt{r}\frac{x}{\ell}\right)
-\frac{1}{2r}$, & in $(\hat{x}_i,L)$.
\end{tabular}
\right.
\end{equation*}
The update 
\begin{equation*}
\alpha_{i+1}(x):=\left\{
\begin{tabular}{ll}
$\displaystyle \frac{1}{4\ell^2}x^2+c_1x+c_2$, & in $(0,\hat{x}_i)$,   \\ [0.2cm]
$\displaystyle b_1\exp\left(-\sqrt{r}\frac{x}{\ell}\right) + b_2\exp\left(\sqrt{r}\frac{x}{\ell}\right)
-\frac{1}{2r}$, & in $(\hat{x}_i,L)$,
\end{tabular}
\right.
\end{equation*}
is explicitly independent on $\alpha_i$ (implicitly, $\alpha_i$ is still present via $\hat{x_i}$). The constants $c_1,c_2$ and $b_1,b_2$ are to be calculated using the boundary conditions in (\ref{bvp_r}), as well as the imposed continuity of both $\alpha_{i+1}$ and $\alpha_{i+1}^\prime$ at $x=\hat{x}_i$.

The ultimate representation for the approximate solution of the boundary value problem (\ref{bvp_r}) obtained at iteration $i\geq0$ is as follows:
\begin{equation}
\alpha_{i+1}(x):=\left\{
\begin{tabular}{ll}
$\displaystyle \frac{1}{4\ell^2}x^2+R(\hat{x}_i)x+1$, & in $(0,\hat{x}_i)$,   \\ [0.2cm]
$\displaystyle T(\hat{x}_i)\exp\left(-\sqrt{r}\frac{x}{\ell}\right) -\frac{1}{2r}$, & in $(\hat{x}_i,L)$,
\end{tabular}
\right.
\label{ai_r}
\end{equation}
where
\begin{equation}
R(\hat{x}_i):=-\frac{ \frac{1}{2\ell^2}\hat{x}_i+\frac{\sqrt{r}}{\ell}\left( 1+\frac{1}{2r}-\frac{1}{4\ell^2}\hat{x}_i^2 \right)}
{1+\frac{\sqrt{r}}{\ell}\hat{x}_i},
\label{Ri}
\end{equation}
and
\begin{equation}
T(\hat{x}_i):=\frac{1+\frac{1}{2r}-\frac{1}{4\ell^2}\hat{x}_i^2}
{1+\frac{\sqrt{r}}{\ell}\hat{x}_i}.
\label{Ti}
\end{equation}
Also, solving $\alpha_{i+1}(x)=0$, one obtains
\begin{equation}
\hat{x}_{i+1}:=\frac{\ell}{\sqrt{r}}
\ln\left( 
\frac{1+2r-\frac{r}{2\ell^2}\hat{x}_i^2}{1+\frac{\sqrt{r}}{\ell}\hat{x}_i} 
\right)
+\hat{x}_i.
\label{xip1}
\end{equation}
which is to be used for constructing the solution at the next iteration. The sketch outlining the construction idea of $\alpha_{i+1}$ in (\ref{ai_r}) using $\alpha_i$ and the corresponding point $\hat{x}_i$ is depicted in Figure \ref{Fig21}(a). The induced point $\hat{x}_{i+1}$ is also depicted in this figure. Note that for initiating the recursion given by (\ref{ai_r}), (\ref{Ri}), (\ref{Ti}) and (\ref{xip1}), one needs $\hat{x}_0$. One of the options that we employ here is to take $\alpha_0(x)=\frac{1}{4\ell^2}x^2-\frac{L}{2\ell^2}x+1$ -- see the left plot in Figure \ref{Fig20}(a) -- such that the required $\hat{x}_0$ reads $\hat{x}_0:=L-\sqrt{L^2-4\ell^2}$. 

For $L:=4\ell$ and $r\in\{4,40\}$ being set in (\ref{ai_r}), Figure \ref{Fig21}(b) depicts the behavior of some first terms of the sequence $\{\alpha_{i+1}\}$, $i\geq0$ with respect to the corresponding reference solution $\alpha^\star$ in (\ref{aStar1d_Lin}). For fixed $r$, there is a numerical evidence that the sequence $\{\alpha_{i+1}\}$ converges. Moreover, the larger $r$, the 'better' the sequence converges to the reference profile $\alpha^\star$, as expected. We now prove the observed trend analytically. 

Indeed, it can be shown that for any $r,\ell>0$ the sequence $\{x_{i+1}\}$, $i\geq0$ defined by (\ref{xip1}) is a monotonically increasing sequence, whose limiting point is
\begin{equation}
\hat{x}^\star:=-\frac{\ell}{\sqrt{r}}+\ell\sqrt{\frac{1}{r}+4}.
\label{xStar}
\end{equation}
It is also obvious that for any fixed $r,\ell>0$, $x^\star\in(0,2\ell)$ and, moreover, it holds that $\lim_{r\rightarrow+\infty}\hat{x}^\star=2\ell$. Inserting (\ref{xStar}) into (\ref{ai_r}), (\ref{Ri}), (\ref{Ti}), we obtain the limit of the sequence $\{\alpha_{i+1}\}$, $i\ge0$, which then represents the (exact!) solution to the boundary value problem (\ref{bvp_r}), that is, the minimizer of (\ref{Es1d_rho}):
\begin{equation*}
\alpha^\star_\rho(x):=\left\{
\begin{tabular}{ll}
$\displaystyle \frac{1}{4\ell^2}x^2+R(\hat{x}^\star)x+1$, & in $(0,\hat{x}^\star)$,   \\ [0.2cm]
$\displaystyle T(\hat{x}^\star)\exp\left(-\sqrt{r}\frac{x}{\ell}\right) -\frac{1}{2r}$, & in $(\hat{x}^\star,L)$,
\end{tabular}
\right.
\end{equation*}
where
\begin{equation*}
R(\hat{x}^\star):=-\frac{1}{2\ell}\sqrt{\frac{1}{r}+4},
\end{equation*}
and
\begin{equation*}
T(\hat{x}^\star):=\frac{1}{2r}\exp\left( -1+\sqrt{1+4r} \right).
\end{equation*}


\end{document}